\documentclass[a4paper]{article}
\usepackage[utf8]{inputenc}
\usepackage[T1]{fontenc}
\usepackage[english]{babel}
\usepackage{enumitem,enumerate}%lists
\usepackage{graphicx,caption,subcaption}%figures,captions and subfigures
\usepackage[retainorgcmds]{IEEEtrantools}
\usepackage{dashrule,xcolor}

\usepackage{mathtools}%math
\def\R{\mathbb{R}}

\def\0{\mathbf{0}}
\def\cross[#1]{\left[ #1 \times \right]}

\usepackage{isomath}% default greek: italic (even for capital), use \mathbf to get upright greek
\usepackage{bm}
\usepackage{amssymb}
\SetMathAlphabet{\mathbf}{normal}{OML}{mdbch}{b}{n}

\def\dbody{\bm{\mathsf{d}}}

\def\ebody{\bm{\mathsf{e}}}

\def\cbody{\bm{\mathsf{c}}}

\def\ubody{\bm{\mathsf{u}}}
\def\xbody{\bm{\mathsf{x}}}
\def\Abody{\bm{\mathsf{A}}}

\def\omegabody{\mathbf{\omega}} % \mathbf{\omega}

\def \abody {\bm{\mathsf{a}}}

\def\fbody{\bm{\mathsf{f}}}

\def\Drag{\mathbb{D}}

\def\dragbody{\mathsf{D}}

\def\Mob{\mathbb{M}}

\def\mobbody{\mathsf{M}}

\def\gbody{\bm{\mathsf{g}}}
\def\xbody{\bm{\mathsf{x}}}

\def\mlab{\bm{m}}
\def\mbody{\bm{\mathsf{m}}}

\def\Bbody{\bm{\mathsf{B}}}

\def\ubody{\bm{\mathsf{u}}}

\def\Pbody{\mathsf{P}}
\def\betabody{\mathbf{\beta}}
\def\etabody{\mathbf{\eta}}
\def\hbody{\bm{\mathsf{h}}}

\def\eps{\varepsilon}

\newcommand{\re}[1]{~\eqref{#1}}
\newcommand{\ret}[2]{~(\ref{#1},\,\ref{#2})}
\newcommand{\e}{\ebody_3}
\newcommand{\ep}[1]{\ebody^{[#1]}}
\newcommand{\bp}[1]{\Bbody^{[#1]}}
\newcommand{\up}[1]{\ubody^{[#1]}}
\newcommand{\gp}[1]{\gbody^{[#1]}}
\newcommand{\taup}[1]{\tau^{[#1]}}

\newcommand{\bma}{\begin{pmatrix}}
\newcommand{\ema}{\end{pmatrix}}

\newcommand{\be}[1]{\begin{equation}\label{#1}}
\newcommand{\ee}{\end{equation}}
\newcommand{\bes}[1]{\begin{equation}\label{#1} \begin{split}}
\newcommand{\ees}{\end{split} \end{equation}}
\newcommand{\ben}[1]{\begin{eqnarray}\label{#1}}
\newcommand{\een}{\end{eqnarray}}

\newcommand{\ue}[1]{\ubody^{[#1]}}

\newcommand{\Rex}[1]{R^{[#1]}}

\definecolor{gray1}{gray}{0.5}
\definecolor{gray2}{gray}{0}
\definecolor{magFrame}{rgb}{0.8667,0,0.5765}

\title{Asymptotic Dynamics of Magnetic Micro-Swimmers}
\author{Pauline Rüegg-Reymond and 
	Thomas Lessinnes}

\begin{document}
\maketitle

%\tableofcontents

\begin{abstract}
Micro-swimmers put into motion by a rotating magnetic field have provided interesting challenges both in engineering and in mathematical modelling. We study here the dynamics of a permanent-magnetic rigid body submitted to a spatially-uniform steadily-rotating magnetic field in Stokes flow. This system depends on two external parameters: the Mason number, which is proportional to the angular speed of the magnetic field and inversely proportional to the magnitude of the field, and the conical angle between the magnetic field and its axis of rotation. This work focuses on asymptotic dynamics in the limits of low and high Mason number, and in the limit of low conical angle. Analytical solutions are provided in these three regimes. In the limit of low Mason number, the dynamical system admits a periodic solution in which the magnetic moment of the swimmer tends to align with the magnetic field. In the limit of large Mason number, the magnetic moment tends to align with the average magnetic field, which is parallel to the axis of rotation. Asymptotic dynamics in the limit of low conical angle allow to bridge these two regimes. Finally, we use numerical methods to compare these analytical predictions with numerical solutions.
\end{abstract}

%----------------------------------------------------------------
\section{Introduction}

Understanding motion at low Reynolds number has been an active research topic for decades~\cite{Brenner1963,Rotne1969,Purcell1977,Happel1983,Rathore2010,Kim2013}. Intuitively, this corresponds to the limit of either very small bodies moving in water, or bodies in a very viscous fluid. It was initially motivated by the study of micro-organisms~\cite{Berg1973,Lighthill1976,Lauga2009}, and more recently also found applications in engineering with the advent of artificial micro-swimmers~\cite{Honda1996,Keaveny2008,Peyer2013a}. These consist in micrometer to millimetre scale devices immersed in fluid, propelled either by a chemical fuel~\cite{Ebbens2010,Gallino2018}, or by an external power source, often an external magnetic field~\cite{Abbott2009,Hwang2012}. An example of application is the precise delivery of a microscopic payload to a specific location in a complicated geometrical environment and through fluids of different rheological properties. In particular, this problem occurs in bio-medical applications, such as targeted drug delivery and microsurgery~\cite{Nelson2010}.

In both~\cite{Ruegg-Reymond2018} and the present study, we focus on the setup where the swimmer is a rigid permanent magnet that experiences a rotating external field, see also~\cite{Man2013,Ghosh2013,Meshkati2014,Fu2015,Morozov2017}. Our theoretical results are valid for rigid swimmers of arbitrary shapes, while our numerical experiments are performed on helices with circular cross-sections.

The swimmer is controlled by adjusting the Mason number $a$ and the angle $\psi$ between the magnetic field and its axis of rotation. The Mason number is a non-dimensional number that represents the balance between the drag and the magnetic load on the swimmer; it will help to keep in mind that it is proportional to the angular velocity of the rotating magnetic field, and inversely proportional to its magnitude. At any moment, the magnetic moment of the swimmer tries to align with the external field but this rotation is slowed by the ensuing Stokes flow. In~\cite{Ruegg-Reymond2018}, this system is modelled as a first-order nonlinear ODE in SO(3) -- see also~\cite{Man2013,Ghosh2013,Morozov2017}. Depending on $a$ and $\psi$ it has exactly zero, one, or two stable relative equilibria. For these motions, the body rotates at the same speed as the external field with its magnetic moment lagging behind the external field. If the body is chiral, this rotation generates a controllable translational motion. Finally, during that study, we observed that the system can also remain out of equilibrium. Experimentally these states would be classified as beyond step-out behaviours.  It is the purpose of the present paper to study them analytically. 

Mathematically, we study this non-linear equation by a combination of analytical and numerical methods. We identify a certain scale $a_c$ such that, in the limit of small Mason numbers $a \ll a_c$, an asymptotic expansion shows that the system can be simplified to a single first-order ODE on the circle: all out of equilibrium solutions are periodic in this limit. In the limit of large Mason numbers $a\gg a_c$, we use the averaging method~\cite{Sanders2007} to obtain an effective equation on SO(3) called the guiding system. Its analysis show that (unsurprisingly) the magnetic moment aligns with the averaged field. But higher order terms coming from the averaging method reveal a secular motion: the swimmer slowly rotates -- timescale of order $a_c/a$ -- about the mean field. To bridge between the low $a$ and large $a$ regimes, we study the limit of small $\sin \psi$ which corresponds to the magnetic field being close to either parallel or antiparallel to its axis of rotation. This reveals a continuous change between the low $a$ regime where the magnetic moment of the swimmer is remaining close to the magnetic field, and the large $a$ regime where the magnetic moment is remaining close to the axis of rotation of the magnetic field, and allows us to characterise the average position of the magnetic moment with respect to the magnetic field, and the amplitude of the excursions of the magnetic moment away from its average position.

Direct numerical integrations for small and large $a$ as well as small $\sin \psi$ corroborate the analytical findings for all particular swimmers we tested. Finally, we use numerical continuation methods~\cite{Kuznetsov2004} to show that these regimes transform into one another as parameters are varied. 

The paper is organised as follows. The dynamical system is described in section~\ref{sec-goveq}. The limits of low Mason number $a$, high Mason number and low  $\sin \psi$ are respectively studied in sections 3, 4 and 5. The latter also describes the transition from low $a$ to high $a$.  Finally, section~\ref{sec-dni} contains numerical integrations of the system and discusses how they compare  with our predictions from the three previous sections. Note also that sections 3, 4 and 5 open with a short summary of the main results therein.

%----------------------------------------------------------------
\section{Governing Equations}
\label{sec-goveq}

We consider a neutrally buoyant rigid body immersed in a viscous fluid filling an infinite three-dimensional space. The medium is permeated by a rotating, spatially homogeneous, magnetic field~$\Bbody$ which we call the external field for short. The body is assumed to be a permanent magnet with magnetic moment~$\mbody$. The shape of the body is taken into account only through its mobility matrix
\be{defM}
\mobbody=\bma \mobbody_{11} & \mobbody_{12} \\ \mobbody_{12}^T &\mobbody_{22}\ema \in \R^{6\times6} \, ,
\ee
which is the inverse of the drag matrix giving the Stokes flow approximation of force and torque due to drag as a linear combination of linear and angular velocities (see~\cite{Happel1983,Kim2013} for instance).

In the Stokes flow limit, the angular velocity $\omegabody$ of such a swimmer is approximated as
\begin{equation}
	\label{ovw}
	\omegabody = \mobbody_{22} \cross[\mbody] \Bbody \, ,
\end{equation}
where after non-dimensionalisation and scaling, the magnetic moment $\mbody$ and magnetic field $\Bbody$ are unit vectors. We leave aside the linear velocity, as the dynamics of the swimmer depend only on its orientation -- see~\cite{Man2013,Meshkati2014} for instance; for a detailed derivation, see~\cite{Ruegg-Reymond2018,Ruegg-Reymond2019}.

\subsection{The Swimmer's Orientation}

We assume that the body is equipped with a right-handed and orthonormal material frame~$D=\{\dbody_1,\, \dbody_2,\, \dbody_3\}$. We also use a fixed, right-handed and orthonormal lab frame~$E=\{\ebody_1,\, \ebody_2,\, \ebody_3\}$. Throughout, we make the slight abuse of notation that involves identifying a vector with its components in the body frame $D$, i.e.~$\mbody$ represents the triplet $(\mbody \cdot \dbody_{1}, \mbody \cdot \dbody_{2}, \mbody \cdot \dbody_{3})$ for instance. In particular, the mobility matrix $\mobbody$ and the magnetic moment $\mbody$ are constant in the body frame. This allows to define the constant matrix
\be{defP}
	\Pbody = \mobbody_{22} \cross[\mbody]
\ee
so that\re{ovw} rewrites
\be{defomega}
	\omegabody(t) = \Pbody \, \Bbody(t) \, ,
\ee
where all time dependences are explicit.

We denote by $R \in$~SO(3) the matrix representing the change of basis between the lab frame and the body frame, i.e. $D = E \, R$. By definition of the angular velocity $\omegabody$, this matrix obeys  the following ode on SO(3)
\be{jumptoR}
	\dot R = R \cross[\omegabody],
\ee
with the notation 
\be{defcross}
[\abody \times] =\bma  0 &-\abody\cdot \dbody_3 & \abody\cdot \dbody_2\\
\abody\cdot \dbody_3 &0 &-\abody\cdot\dbody_1\\
-\abody\cdot\dbody_2 &\abody\cdot \dbody_1 &0 \ema. 
\ee

Taking the transposed of\re{jumptoR} and working column by column yields
\be{jumptoes}
	\dot{\ebody}_{i} = - \omegabody \times \ebody_{i} ,
\ee
where as per our convention, the vectors $\ebody_{i}$ of the lab frame also represent the triplets $(\ebody_{i} \cdot \dbody_{j})_{j = 1,2,3}$.

We define a third frame, the magnetic frame $\tilde E$, in which the magnetic field $\Bbody$ and its axis of rotation are fixed. The axis of rotation is also fixed in the lab frame, and without loss of generality we pick it as $\ebody_{3}$, so that
\be{defEtilde}
\widetilde E= E \, R_3(a \, t) \qquad \textrm{ where}\qquad
R_3(\varphi) := \bma \cos \varphi & -\sin \varphi& 0 \\ \sin\varphi &\cos \varphi &0 \\ 0& 0&1 \ema \, ,
\ee
and where $a$ is the Mason number obtained after non-dimensionalisation of the angular speed $\alpha$ of the magnetic field:
\be{defa} 
	a = \frac{\alpha\, \eta\, \ell^3}{m\, B} \, ,
\ee
where $\eta$ is the dynamic viscosity of the fluid, $\ell$ is a characteristic length, $m$ is the magnitude of the magnetic moment, and $B$ is the magnitude of the magnetic field (see~\cite{Ruegg-Reymond2018,Ruegg-Reymond2019} for details).

Finally, we define $Q$ as the matrix representing the change of basis between the magnetic frame and the body frame, i.e. $\tilde E = D \, Q$, which implies
\be{defQ}
	Q(t) =R^T(t)\, R_3(a\, t) \, .
\ee
Thus the magnetic field $\Bbody$ can be written as
\be{defB}
	\Bbody(t) = Q(t) \, \begin{bsmallmatrix} \sin \psi \\ 0 \\ \cos \psi \end{bsmallmatrix} \, ,
\ee
where $\psi > 0$ is the conical angle between the magnetic field and its axis of rotation -- note that this implies that the lab frame is chosen so that the magnetic field lies in the $(\ebody_{1},\ebody_{3})$-plane and that $\ebody_{1} \cdot \Bbody > 0$ at time $t = 0$, which can be assumed without loss of generality. In this picture, taking a time derivative of\re{defQ} and using\re{jumptoR}, yields
\be{jumptoQ}
	\dot Q=  \left[\left(a \, Q\, \begin{bsmallmatrix} 0 \\ 0 \\1 \end{bsmallmatrix} - \omegabody\right)\times\right]\, Q \, .
\ee

It is also straightforward to show that if we know $\ebody_3$ and $\Bbody$, then $\ebody_{1}$, $\ebody_{2}$, and $\ebody_{3}$ are known (provided $\sin \psi \neq 0$). In this picture, it is advantageous to use\re{defB} to replace\re{jumptoes} by
\be{jumptoeB}
\left\{\begin{split}
\dot \ebody_3 &= - \omegabody \times \ebody_3 \, , \\
\dot \Bbody &= (a \, \ebody_3 - \omegabody) \times \Bbody \, .
\end{split}\right.
\ee

The three systems of equations\re{jumptoes},\re{jumptoQ} and\re{jumptoeB} are equivalent. Using\re{defomega} and\re{defB}, they can be written in closed form respectively as
\begin{align}
	&\dot{\ebody}_{i} = - (\Pbody \, \Bbody) \times \ebody_{i} \, ,\text{ where }  \Bbody = \sin \psi \, (\cos(a \, t) \, \ebody_{1} + \sin (a \, t) \, \ebody_{2}) + \cos \psi \, \ebody_{3} \, ,
	\label{syses}
	\end{align}
and where the $\ebody_i$ form a right-handed orthonormal basis;
\begin{align}
	&\dot Q=  \left[\left(a\, Q \,\begin{bsmallmatrix} 0  \\ 0 \\ 1 \end{bsmallmatrix} - \Pbody Q \begin{bsmallmatrix} \sin \psi \\ 0 \\ \cos \psi \end{bsmallmatrix}\right)\times\right]\, Q \, ,
	\label{sysQ}
	\end{align} 
	where $Q\in\textrm{SO}(3)$;
	\begin{align}
	&\begin{cases}
	\dot \ebody_3 = - (\Pbody \Bbody) \times \ebody_3 \, , \\
	\dot \Bbody = (a\, \ebody_3 - \Pbody \, \Bbody) \times \Bbody \, ,
	\end{cases}
	\label{syseB}
\end{align}
where $\ebody_3\cdot \Bbody=\cos\psi,\, \ebody_3\cdot \ebody_3=1$ and\footnote{This last equation comes from the scaling already mentioned under equation~(\ref{ovw}).} $\Bbody\cdot \Bbody=1$.

Both the systems\re{syseB} and\re{sysQ} are autonomous, whereas in\re{syses} $\Bbody$ depends explicitly on $t$. 
When solving analytically, we found it easiest to work with $\ebody_3$ and $\Bbody$ in section~\ref{sec-lowa}, and with the $\ebody_{i}$s in sections~\ref{sec-higha} and~\ref{sec-smallpsi}. When solving numerically in section~\ref{sec-dni}, we found it easiest to work with $Q$.

The behaviour of the swimmer therefore depends on the two parameters $a$ and $\psi$. The aim of this paper is to study its asymptotic dynamics in the limits of $a \to 0$, $a \to \infty$, and $\sin\psi \to 0$. The steady states of these systems are studied in detail in~\cite{Ruegg-Reymond2018,Ruegg-Reymond2019}.

\subsection{The $\Pbody$ matrix}\label{sec-Pmat}

We will use the matrix $\Pbody$ defined in\re{defP} and its singular value decomposition throughout our analysis. We introduce it here.

Let $\sigma_{0} = 0$, $\sigma_{1}$, and $\sigma_{2}$ be the singular values of $\Pbody$ with corresponding right-singular vectors $\betabody_{0} = \mbody$, $\betabody_{1}$, and $\betabody_{2}$ and left-singular vectors $\etabody_{0} = \mobbody_{22}^{-1} \mbody / \left\| \mobbody_{22}^{-1} \mbody \right\|$, $\etabody_{1}$, and $\etabody_{2}$. Note that $\mobbody_{22}$ is a symmetric and positive definite block of the mobility matrix\re{defM} so that its inverse exists~\cite{Happel1983}. Furthermore, these definitions imply that $\betabody_0\cdot \etabody_0 >0$ so that there exist angles $\iota\in[0,\pi/2)$ and $\zeta\in[0,2\pi)$ such that 
\begin{equation*}
\etabody_0 = \cos \iota \, \betabody_{0} + \sin \iota ~( \sin \zeta \, \betabody_1 -\cos \zeta\, \betabody_2) \, .
\end{equation*}
We will see in sections~\ref{sec-lowa} and~\ref{sec-higha} that $\iota$ is a crucial material parameter for the existence of out of equilibrium solutions of the system. 

By definition, we have the relations
\begin{align*}
	\Pbody \, \betabody_{i} &= \sigma_{i} \, \etabody_{i} \, , & \Pbody^T \, \etabody_{i} &= \sigma_{i} \, \betabody_{i} \, ,\label{singvecttheory}
\end{align*}
and each set of singular vectors forms an orthonormal basis that we can assume to be right-handed. These bases are also constant in the body frame.

% ----------------------------------------------------------------
\section{Periodic Orbits at Small $a$}\label{sec-lowa}

Low values of $a$ can be achieved experimentally either by considering strong magnetic effects, slowly rotating external fields,  or particularly small bodies -- see\re{defa}. We show that in this case, the magnetic field and the magnetic moment of the swimmer align on a time scale of order 1. Thereafter, the leading order dynamic is completely specified by the first order differential equation\re{dTlambda} on $\mathbb R$. Specifically, we show that \begin{itemize} 
\item~~if $\psi \in (\pi/2-\iota, \pi/2 +\iota)$, then there exists a unique stable relative equilibrium;
\item~~if $\psi \notin [\pi/2-\iota, \pi/2 +\iota]$, then the leading order dynamics exhibits a single periodic orbit of period 
$$
\frac {2\pi \, \cos \iota}{a \, \sqrt{\cos (\iota+ \psi) \, \cos(\iota -\psi) }} \, .
$$
\end{itemize}

To obtain this result, we first note that when $a\ll \mathrm{min}\{1,\sigma_1,\sigma_2\}$, it takes an $\mathcal O(1/a)$ time for the magnetic field to make a full rotation around its axis. Since we are interested in the behaviour of the system after many revolutions, we rescale to a longer time scale $T = a\, t$. The system\re{syseB} becomes
\begin{equation} \label{systemlowarescaled}
\left \{ 
\begin{split} 
a\, \frac {d\ebody_3}{dT} & = - (\Pbody \, \Bbody) \times \ebody_3 \, , \\
a\, \frac {d\Bbody}{dT} & = a \, \ebody_3 \times \Bbody -  (\Pbody \, \Bbody) \times \Bbody \, .
\end{split}
\right .
\end{equation}

In the following, we perform a singular expansion analysis of\re{systemlowarescaled}. In the inner layer, that is when $T=\mathcal O(1/a)$ or equivalently when $t= \mathcal O(1)$, the governing equation are\re{syseB}.
Expanding
\begin{align} \label{expandlowa}
\ebody_3 &=  \ebody_3^{[0]} + a\, \ebody_3^{[1]} + \mathcal O(a^2) \, , & 
\Bbody & = \Bbody^{[0]} + a\, \Bbody^{[1]} + \mathcal O(a^2) \, ,
\end{align}
 substituting\re{expandlowa} in\re{syseB} and matching at each order in $a$ we find that the equation for $\bp 0$ decouples and reads
\begin{equation*}
\dot{\Bbody}^{[0]} = - (\Pbody \,  \bp 0) \times \bp 0.
\end{equation*}
This equation has two equilibria $\bp 0= \pm \betabody_0$ where $\bp 0 = \betabody_0$ is stable and globally attracting while $\bp 0 = -\betabody_0$ is unstable.

In the outer layer, that is when $T>1$, we substitute an expansion of the form\re{expandlowa} in\re{systemlowarescaled} and match at zeroth and first order in $a$ to find 
\begin{equation} \label{sysla}
\left \{ 
\begin {split} 
&(\Pbody \, \bp 0 ) \times \ep{0}_{3}  = 0 \, , \\
&(\Pbody \, \bp 0 ) \times \bp 0  = 0 \, , \\
&\frac {d \ep{0}_{3} } {dT} = - ( \Pbody \bp 1 ) \times \ep{0}_{3}  - ( \Pbody \bp 0 ) \times \ep{1}_{3} \, ,\\
&\frac {d \bp 0 } {dT} = ( \ep{0}_{3} -\Pbody \bp 1 ) \times \bp 0  - ( \Pbody \bp 0 ) \times \bp 1 \, .
\end{split}
\right.
\end{equation}

Because $\Bbody$ and $\e$ are unit vectors and $\Bbody \cdot \e = \cos \psi$ is constant, we have
\be{condlowa}
\begin{aligned} 
	\bp 0 \cdot \bp 0 &= 1 \, , &
	\ep{0}_{3} \cdot \ep{0}_{3} &=1 \, ,  & 
	\bp 1 \cdot \bp 0 &= 0 \, , \\
	\ep{1}_{3} \cdot \ep{0}_{3} &= 0 \, , &
	\ep{0}_{3} \cdot \bp 0 &= \cos \psi \, .
\end{aligned}
\ee

The first two equations in\re{sysla} imply that $\Pbody \, \bp 0 = 0$ and hence 
\begin{equation}\label{bet0lowa}
\bp 0 = \pm \betabody_0 \, ,
\end{equation} 
where we choose the `$+$' sign because it corresponds to the attracting equilibrium of the inner layer. Substituting this result in the last equation of\re{sysla} yields 
\begin{equation}\label{lowafourthagain}
\Pbody \, \bp 1 \times \bp 0 = \ep{0}_{3} \times \bp 0 \, .
\end{equation}

Furthermore, because of\ret{condlowa}{bet0lowa}, there exists a function $\lambda(T)$ such that 
\begin{equation}\label{deflambdalowa}
\ep{0}_{3} (T) = \cos \psi\, \betabody_0 + \sin \psi \, \big(\cos \lambda(T) \, \betabody_1 +\sin \lambda(T) \, \betabody_2\big).
\end {equation}
 Next, equations~(\ref{condlowa}-\ref{deflambdalowa}) together with $\Pbody \, \bp 1\perp \etabody_0$ yield
\begin{equation}\label{solvepb1lowa}
\Pbody \, \bp 1(T)=\sin \psi \, \Big( \cos \lambda(T) \, \betabody_1 + \sin \lambda(T) \, \betabody_2 + \sin \big( \lambda(T) - \zeta\big) \tan \iota \, \betabody_0 \Big ) \, ,
\end{equation}
where $\iota \in \left[ 0, \pi/2 \right)$ and $\zeta\in[0,2\pi]$ have been defined in section~\ref{sec-Pmat}.

Finally, substituting\ret{deflambdalowa}{solvepb1lowa} in the third equation of\re{sysla} yields 
\begin{equation}\label{dTlambda}
\frac {d \lambda}{dT} = \cos \psi - \sin \psi \, \tan \iota~ \sin (\lambda - \zeta) \, .
\end{equation}

When $\pi/2-\iota < \psi < \pi/2+\iota$, we have $|\cos \psi|<\sin \psi \, \tan \iota$ and\re{dTlambda} has two equilibria given by
\begin{equation*}
	\lambda= \zeta +  \frac \pi 2  \pm \arccos \left( \frac{\cos \psi \cos \iota}{\sin \psi \sin \iota} \right) \, .
\end{equation*}
The equilibrium with a minus sign in $\pm$ is stable and the other one is unstable. 

If $|\psi -\pi/2|> \iota$,  then $|\cos \psi|>\sin \psi \, \tan \iota$ and $d\lambda / dT$ never changes sign. In consequence, the leading order dynamic exhibits a periodic solution of period 
\begin{equation}
\label{periodlowa}
\int_{\zeta-\pi}^{\zeta + \pi} \frac {d\lambda} {\cos \psi - \sin \psi \, \tan \iota \, \sin (\lambda - \zeta)} = \frac {2 \pi \, \cos \iota}{\sqrt{\cos (\iota+ \psi) \, \cos(\iota -\psi) }} \, .
\end{equation}
In the body frame, this periodic dynamic corresponds to the axis of rotation $\ebody_{3}$ of the magnetic field itself rotating about $\betabody_0$; clockwise when $\psi\in(0,\pi/2-\iota)$ and anti-clockwise when $\psi\in(\pi/2+\iota , \pi)$. 
This bifurcation between stable equilibria for $\psi \in [\iota, \pi-\iota]$ and periodic orbit when $\psi \notin [\iota, \pi-\iota]$ occurs through a periodic solution of infinite period.

%----------------------------------------------------------------
\section{Asymptotic Dynamics at Large $a$}\label{sec-higha}

Large values of the Mason number $a$ correspond in experiments to either weak or rapidly rotating magnetic fields. We show that in this setting, the magnetic moment $\mbody$ tends to align with the average magnetic field, that corresponds either to its axis of rotation $+\ebody_{3}$ or to the opposite of its axis of rotation $-\ebody_{3}$ depending on the conical angle $\psi$. The mismatch between the magnetic moment and $\pm \ebody_{3}$ is of order $1/a$. We also show that there is a slow residual rotation of the swimmer about the average field, with period of order $a$.

To analyse the case of $a\gg1$, we work with the version\re{syses} of our system
\begin{equation} \label{govhigha}
\dot{\ebody}_i  = - (\Pbody \, \Bbody) \times \ebody_i \, , 
\end{equation}
and remember that the external field is given as rotating at constant angular velocity with respect to the lab frame so that
\begin{equation} \label{externalfield}
\Bbody(t)= \sin \psi\, \big(\cos (a\, t)\, \ebody_1(t) + \sin (a\, t)\, \ebody_2(t) \big) + \cos \psi\, \ebody_3(t).
\end{equation} 

We will analyse\re{govhigha} by applying the averaging method described in~\cite[]{Sanders2007}. The main idea is that because $\Bbody$ changes much faster than $\ebody_3$, we can approximate the effect of $\Bbody$ on the dynamic of the system by averaging $\Bbody$ over one of its period of revolution $\frac{ 2\pi}{a}$. The method transforms the non-autonomous system\ret{govhigha}{externalfield} into an autonomous averaged differential equation called the guiding system. The averaging procedure is carried out in Section~\ref{app-avmeth} and the resulting guiding system is studied in Section~\ref{app-guide}.

\subsection{Averaged Governing Equations} \label{app-avmeth}

We first rescale time to $T = a\, t$ and define $\varepsilon = 1/a\ll1$ so that\re{govhigha} becomes 
\begin{equation} \label{govhighascaled}
\frac {d \ebody_i } {dT} = -\varepsilon \,  (\Pbody \,  \Bbody) \times \ebody_i \, . 
\end{equation}

The averaging operator, noted by an overline, is defined as follows: to any function $X(\ebody_1,\, \ebody_2,\, \ebody_3,\,T)$ $2\pi$-periodic in $T$, it associates the averaged function
\begin{equation*} 
\overline{X}(\ebody_1,\, \ebody_2,\, \ebody_3) = \frac 1 {2 \pi} \int_0^{{2\pi}} X(\ebody_1,\, \ebody_2,\, \ebody_3, T) \, dT \, ,
\end{equation*}
where the integration is performed while keeping the $\ebody_i$s constant. We also define the function $\Bbody^\star$ which depends on three vectors and a scalar
\begin{equation*}
\Bbody^\star(\xbody_1,\xbody_2,\xbody_3,T) = \sin \psi\, \big(\cos (T)\, \xbody_1+ \sin (T)\, \xbody_2 \big) + \cos \psi\, \xbody_3 \, ,
\end{equation*}
such that $\Bbody(T) = \Bbody^\star\big(\ebody_1(T),\ebody_2(T),\ebody_3(T),T\big)$. Note that $\overline {\Bbody^\star}(\xbody_1,\xbody_2,\xbody_3) = \cos \psi \, \xbody_3$.

Given two functions $\up 1$ and $\up 2$, the near identity transformation (see appendix~\ref{app-exp SO3})
\begin{equation}\label{nearidentity}
\ebody_i = \cbody_i + \varepsilon \, \left(\up 1 \times \cbody_i\right)+ \varepsilon^2 \, \left(\up 2 \times \cbody_i +\frac 1 2\, \up 1 \times (\up 1 \times \cbody_i)\right) + \mathcal O(\varepsilon^3)
\end{equation}
transforms the differential equation\re{govhighascaled} into a differential equation for the right-handed orthonormal frame of the $\cbody_i$s: 
\begin{equation}\label{guidesys0}
\frac {d \cbody_i}{dT} = \varepsilon \, \gp 1 \times \cbody_i + \varepsilon^2 \, \gp 2 \times \cbody_i + \mathcal O(\varepsilon^3),
\end{equation} 
where the $\boldsymbol g^{[i]}$ are computed hereunder.

The functions $\up i$ ($i=1,\,2$) appearing in\re{nearidentity} are to be understood as explicit functions of the vectors $\cbody_j$ ($j=1,\,2,\,3$) and of time $T$. We further require that they are $2\pi$-periodic in $T$. Then they can be chosen such that the functions $\gp i$ ($i=1,\,2$) appearing in\re{guidesys0} are independent of time. In that case, solutions of the truncated equation 
\begin{equation} \label{guideT}
\frac {d \cbody_i}{dT} = \varepsilon \, \gp 1 \times \cbody_i + \varepsilon^2 \, \gp 2 \times \cbody_i
\end{equation} 
are guaranteed to remain close to the solutions of\re{govhighascaled} up to an order $\mathcal O(\varepsilon)$ on a timescale of order $T=\mathcal O(1/\varepsilon^2)$~\cite{Sanders2007}. This means that in the original time scale, we are guaranteed an approximation of order $\mathcal O(1/a)$ on a time scale of order $t=\mathcal O (a)$ (remember that $a\gg1$).

Substituting\re{nearidentity} in\re{govhighascaled}, expanding in $\varepsilon$, and matching at each order gives 
\begin{align}
&\left .\partial_T \up 1\right|_{(\cbody_1, \cbody_2, \, \cbody_3,\, T)} =- \Pbody\, \Bbody^\star(\cbody_1, \cbody_2, \, \cbody_3,\, T) - \gp 1 (\cbody_1, \cbody_2, \, \cbody_3) \, ,\label{highaexp}\\
&\partial_T \up 2 =
 - \Pbody \left(\up 1 \times \Bbody^\star \right) 
- \frac{\partial \up 1}{\partial \cbody_j} \cdot \left(\gp 1 \times \cbody_j\right)
- \frac 1 2 \up 1 \times \partial_T \up 1  + \gp 1 \times \up 1 - \gp 2 \, ,~~\label{highaexp2}
\end{align}
with implied summation on repeated indices and where the functions $\up 1$, $\up 2$ and $\Bbody^\star$ appearing in\re{highaexp2} are evaluated at $(\cbody_1, \cbody_2, \, \cbody_3,\, T)$ and the functions $\gp 1$ and $\gp 2$ are evaluated at $(\cbody_1, \cbody_2, \, \cbody_3)$.

Requiring the function $\up{1}$ to be periodic implies that the left-hand side of\re{highaexp} vanishes upon averaging over time $T$. Accordingly, we find that 
\begin{equation} 
\gp 1(\cbody_1, \cbody_2, \, \cbody_3) = - \overline{\Pbody\, \Bbody^\star}(\cbody_1, \cbody_2, \, \cbody_3)  = - \cos \psi\, \Pbody\, \cbody_3 \, .\label{g1}
\end{equation}

After substituting\re{g1} in\re{highaexp}, we integrate and find that
\begin{equation}\label{u1}
\up 1 (\cbody_1, \cbody_2, \, \cbody_3, T) = \sin \psi \, \big (- \sin T\, \Pbody \, \cbody_1 + \cos T\, \Pbody \, \cbody_2) + \Abody(\cbody_1, \cbody_2, \, \cbody_3) \, ,
\end{equation}
where $\Abody$ is an arbitrary function. Here, we choose $\Abody = 0$ so that $\overline {\up 1} =0$.

Substituting\ret{g1}{u1} in\re{highaexp2}, and averaging gives 
\begin{equation}\label{g2}
\gp 2 (\cbody_1, \cbody_2, \, \cbody_3) = \frac {\sin^2 \psi}{2} \Big(\Pbody \, (\cbody_1 \times \Pbody \, \cbody_2)  + \Pbody \, (\Pbody \, \cbody_1 \times \cbody_2) - (\Pbody \, \cbody_1) \times (\Pbody \, \cbody_2)\Big) \, .
\end{equation}

Coming back to the time scale $t=T/a$, and substituting\ret{g1}{g2} in\re{guideT}, the guiding system becomes 
\begin{equation}
	\label{guide}
	\begin{aligned}
		\frac{d \cbody_i}{dt} =& -\cos \psi \, (\Pbody \, \cbody_3)\times \cbody_i \\
		& + \varepsilon \, \frac{\sin^2\psi}{2} \, \Big(\Pbody \, (\cbody_1 \times \Pbody \, \cbody_2)  + \Pbody \, (\Pbody \, \cbody_1 \times \cbody_2) - (\Pbody \, \cbody_1) \times (\Pbody \, \cbody_2)\Big)\times \cbody_i \, .
	\end{aligned}
\end{equation}

\subsection{Analysis of the Guiding System}\label{app-guide}
The argument of the previous section transformed the non-autonomous system\ret{govhigha}{externalfield} into the approximate autonomous system\re{guide}. In this section we show that the first-order solution of\re{guide} always exhibits stable periodic motion. 

First, remark that if we truncate the guiding system\re{guide} to zeroth order in $\varepsilon$, we obtain 
\begin{equation}\label{guide0}
\frac{d \cbody_i}{dt} = -\cos \psi \, (\Pbody \, \cbody_3)\times \cbody_i \, .
\end{equation}
The system\re{guide0} has two families of one-parameter equilibria of the form\footnote{Recall that according to section~\ref{sec-Pmat} $\Pbody\, \cbody_3=0$ iff $\cbody_3= \pm \betabody_0$} 
\begin{equation} \label{zerostable}
\left \{
\begin{split}
\cbody_1 &= \cos \chi \, \betabody_1 + \sin \chi \, \betabody_2 \, ,\\
\cbody_2 &=  -\sin \chi\, \betabody_1 + \cos \chi\, \betabody_2 \, ,\\ 
\cbody_3 & =  \pm \betabody_0 \, .
\end{split} 
\right .
\end{equation}

It is straightforward to show that all equilibria corresponding to $\cbody_3= \textrm{Sign}(\cos \psi) \, \betabody_0$ are meta-stable in the sense that all eigenvalues of the associated stability matrix $|\cos \psi| \, \Pbody \, \cross[\betabody_{0}]$ are strictly negative but for one that vanishes and corresponds to motion within the continuous family. It is also straightforward to show that the stable manifold is almost globally attracting -- that is it is attracting for all initial values that do not lie strictly on the manifold of unstable equilibria.

However as the system gets near this zeroth-order stable manifold, the magnitude of the zeroth order term in\re{guide} decreases, and the first order term can no longer be neglected. For the long term behaviour of the system, we therefore expect the system to be close to but not quite on the equilibrium\re{zerostable} with the sign chosen so as to match that of $\cos \psi$. We therefore define a function 
\begin{equation*}
\tau(t) = \tau^{[0]} (t) + \varepsilon\, \tau^{[1]}(t) + \mathcal O(\varepsilon^2),
\end{equation*}
 that specifies elements in the stable manifold via
\begin{equation}
	\label{fasbetas}
	\left\{\begin{split}
		\fbody_1(t) &= \cos \tau(t) \, \betabody_1 + \sin\tau(t) \, \betabody_2 \, ,\\
		\fbody_2(t) &=  \varsigma\,\big(-\sin \tau(t)\, \betabody_1 + \cos \tau(t)\, \betabody_2 \big) \, ,\\ 
		\fbody_3(t) & =  \varsigma\, \betabody_0 \, ,
	\end{split}\right.
\end{equation} 
where $\varsigma=\textrm{Sign}(\cos \psi)$. We then expand the $\cbody_i$s as follows
\begin{equation}
	\label{expandc}
	\cbody_i = \fbody_i + \varepsilon \, \xbody \times \fbody_{i} + \mathcal O (\varepsilon^2) \, ,
\end{equation} 
where $\xbody(t) = x_1(t) \,\fbody_1(t) + x_2(t) \, \fbody_2(t)$ is to be determined by the differential equation\re{guide}.

Substituting\re{expandc} in\re{guide} leads to
\begin{equation}
	\label{expandguide}
	\begin{aligned}
	&\begin{aligned}
		\varsigma \, \dot \tau^{[0]} \, \fbody_3\times \fbody_i +  \varepsilon \, \Big( & \big( \varsigma \, \dot\tau^{[1]} \, \fbody_3  + \dot x_1 \, \fbody_1 + \dot{x}_{2} \, \fbody_{2} \big) \times \fbody_i \\
		& + \, \varsigma\, \dot \tau^{[0]} \, \big( (\fbody_3 \times \xbody) \times \fbody_i +  \xbody \times(\fbody_3\times \fbody_i) \big) \Big) + \mathcal O(\varepsilon^2) \\
	\end{aligned}\\
	&\begin{aligned}
		\qquad = \varepsilon \, \Big( - \cos \psi \, \Pbody \, (\xbody \times \fbody_{3}) + \gp 2 (\fbody_1, \fbody_2, \, \fbody_3) \Big) \times \fbody_{i} + \mathcal O(\varepsilon^2) \, ,
	\end{aligned}
	\end{aligned}
\end{equation}
where $\gp 2 (\fbody_1, \fbody_2, \, \fbody_3)$ is given by\re{g2} where substitution according to\re{fasbetas} yields 
\begin{equation*}
	\gp 2 (\fbody_{1} , \fbody_{2}, \fbody_{3}) = \varsigma \, \frac{\sin^2 \psi}{2} \, \big( \Pbody \, (\betabody_{1} \times \Pbody \, \betabody_{2} - \betabody_{2} \times \Pbody \, \betabody_{1} ) - \sigma_{1} \, \sigma_{2} \, \etabody_{0} \big) \, ,
\end{equation*}
which is independent from $\tau$.

Matching orders in\re{expandguide}, we find at zeroth order  $\dot \tau^{[0]}= 0$ and at first order
\begin{equation}
	\label{govguideo1}
	\left\{\begin{split}
		\varsigma \, \dot{\tau}^{[1]} &= \hbody (\fbody_{1}) \cdot \fbody_{2} \\
		\dot{x}_{1} &= \hbody (\fbody_{2}) \cdot \fbody_{3} \\
		\dot{x}_{2} &= \hbody (\fbody_{3}) \cdot \fbody_{1} \, ,
	\end{split}\right.
\end{equation}
where 
\begin{equation*}
	\hbody (\fbody_{i}) = \Big( - \cos \psi \, \Pbody \, (\xbody \times \fbody_{3}) + \gp 2 (\fbody_1, \fbody_2, \, \fbody_3) \Big) \times \fbody_{i} \, .
\end{equation*}
Because $\dot \tau^{[0]} =0$, the function $\fbody_i(t)$ evolves on a slow time scale $\mathcal O(1/\varepsilon)$. Accordingly we  look for equilibria of the equations in $(x_1,\, x_2)$ while keeping all $\fbody_{i}$s constant. We find a single equilibrium of the form
\begin{equation*}
	\begin{bmatrix} x_{1} \\ x_{2} \end{bmatrix} =
	\frac{1}{\sigma_{1} \, \sigma_{2} \, \cos \psi \, \cos \iota}
	\begin{bmatrix}
		\Pbody \, \fbody_{1} \cdot \big( \fbody_{3} \times \gp 2 (\fbody_{1}, \fbody_{2} , \fbody_{3}) \big) \\
		\Pbody \, \fbody_{2} \cdot \big( \fbody_{3} \times \gp 2 (\fbody_{1}, \fbody_{2} , \fbody_{3}) \big)
	\end{bmatrix}
\end{equation*}
which leads to
\begin{equation*}
	\xbody = \frac{1}{\sigma_{1} \, \sigma_{2} \, \cos \psi \, \cos \iota} \, ( \mathbb{I} - \fbody_{3} \, \fbody_{3}^T) \, \Pbody^T \, \big(\fbody_{3} \times \gp 2 (\fbody_{1} , \fbody_{2}, \fbody_{3}) \big) \, .
\end{equation*}

Substituting $\xbody$ accordingly, the first equation of\re{govguideo1}  becomes
\begin{equation*}
	\dot{\tau}^{[1]} = - \varsigma \, \sigma_{1} \, \sigma_{2} \, \frac{\sin^2 \psi}{2 \cos \iota} \, ,
\end{equation*}
so that
\begin{equation}
	\label{tau eq}
	\tau(t) = \tau_{0} - \eps \, \varsigma \, \sigma_{1} \, \sigma_{2} \, \frac{\sin^2 \psi}{2 \cos \iota} \, t + \mathcal O (\eps^2) \, .
\end{equation}

In conclusion, we have shown that in the limit of large Mason number $a$, the magnetic moment $\mbody$ tends to align with the average magnetic field, which is $\pm\ebody_{3}$ depending on the sign of $\cos \psi$. The mismatch between $\mbody$ and $\pm \ebody_{3}$ is of order $\eps$. Indeed gathering our findings\footnote{Remember that $\betabody_{0} = \mbody$.} we have
\begin{equation*}
	\ebody_{3} = \varsigma \, \mbody + \mathcal O (\eps) \, ,
\end{equation*}
where $\varsigma = \mathrm{sign} \, \cos \psi$. Furthermore,
\begin{align*}
	\ebody_{1} &= \cos \tau(t) \, \betabody_{1} + \sin \tau(t) \, \betabody_{2}  + \mathcal O (\eps) \\
	\ebody_{2} &= - \varsigma \,  \sin \tau(t) \, \betabody_{1} + \varsigma \, \cos \tau(t) \, \betabody_{2} + \mathcal O (\eps) \, ,
\end{align*}
where $\tau(t)$ is given by\re{tau eq}, so that viewed from the lab frame, there is a slow residual rotation of the body frame about the average field.

%----------------------------------------------------------------
\section{Asymptotic Dynamics at Small $\sin\psi$}\label{sec-smallpsi}

In this section, we analyse the case of a magnetic field $\Bbody$ almost parallel to its axis of rotation $\ebody_{3}$; this corresponds to $\sin \psi \ll 1$. We show that in the small $\sin \psi$ regime, the magnetic moment describes a circle in the magnetic frame, whose centre shifts from the time-dependent magnetic field $\Bbody$ to the average magnetic field $\pm \ebody_{3}$ as the Mason number $a$ goes from asymptotically small to asymptotically large, and whose radius goes to zero in both limits $a \to 0$ and $a \to \infty$. This regime thus bridges the small $a$ regime studied in section~\ref{sec-lowa} and the large $a$ regime studied in~\ref{sec-higha}.

We analyse the case $\sin\psi \ll 1$ by setting $\eps = \sin \psi$ and performing an asymptotic expansion in $\eps$. Note that we have $\sin \psi \ll 1$ both for $\psi \ll 1$ and for $\pi-\psi \ll 1$, that is for $\Bbody$ close to either of $\pm \ebody_{3}$. The equation for the dynamics of the lab frame\re{syses} becomes
\be{sys-smallpsi}
	\dot{\ebody}_i  = - \varsigma \, (\Pbody \, \ebody_{3}) \times \ebody_i
	- \eps \, (\Pbody \, R_{3}(a \, t) \, \ebody_{1}) \times \ebody_i
	- \varsigma \, \frac{\eps^2}{2} \, (\Pbody \, \ebody_{3}) \times \ebody_i + \mathcal O(\eps^3) \, ,
\ee
where $\varsigma = \mathrm{sign} (\cos \psi)$.

\subsection{Asymptotic expansion}\label{sec-smallpsi exp}
\subsubsection{Zeroth order}

The zeroth order dynamics is given by the equation
\begin{equation*}
	\dot{\ebody}_i  = - \varsigma \, ( \Pbody \, \ebody_{3} ) \times \ebody_i \quad \text{ for } i = 1,2,3 \, ,
\end{equation*}
which is in equilibrium for $\Pbody \, \ebody_{3} = 0$, i.e. $\ebody_{3} = \pm \betabody_{0}$. This implies that we have two families of equilibria given by
\begin{align}
	\label{zeroth-eq-smallpsi}
	\ebody_{1} &= \cos \tau \, \betabody_{1} \pm \sin \tau \, \betabody_{2} \, , &
	\ebody_{2} &= -\sin \tau \, \betabody_{1} \pm \cos \tau \, \betabody_{2} \, , &
	\ebody_{3} &= \pm \betabody_{0} \, ,
\end{align}
where $\tau$ is a parameter. Furthermore, the equilibria for which the $\pm$ sign is $\varsigma$ are stable. As in section~\ref{app-guide}, the magnitude of the zeroth order term in\re{sys-smallpsi} decreases as the system approaches this zeroth-order stable manifold. To study the long-term behaviour of the system, we must therefore also take into account higher order terms, as we expect the solutions to\re{sys-smallpsi} to be close to but not quite on the equilibrium\re{zeroth-eq-smallpsi} with the sign matching $\varsigma$. To this end we define a function
\begin{equation*}
	\tau(t) = \taup{0}(t) + \eps \, \taup{1}(t) + \eps^2 \, \taup{2}(t) + \mathcal O (\eps^3) \, ,
\end{equation*}
that specifies elements of the stable manifolds as
\begin{equation*}
	\left\{\begin{split}
		\ep{0}_{1}(t) &= \cos \tau (t) \, \betabody_{1} + \varsigma \, \sin \tau (t) \, \betabody_{2} \\ 
		\ep{0}_{2}(t) &= -\sin \tau (t) \, \betabody_{1} + \varsigma \, \cos \tau (t) \, \betabody_{2} \\
		\ep{0}_{3}(t) &= \varsigma \, \betabody_{0} \, .
	\end{split}\right.
\end{equation*}
We then expand the $\ebody_{i}$s as
\begin{equation}
	\label{smallpsi exp es}
	\ebody_{i} = \ep{0}_{i} + \eps \, \ep{1}_{i} + \eps^2 \, \ep{2}_{i} + \mathcal O (\eps^3) \, ,
\end{equation}
where
\begin{align*}
	\ep{1}_{i} &= \up{1} \times \ep{0}_{i} \, , & \ep{2}_{i} &= \up{2} \times \ep{0}_{i} + \frac{1}{2} \, \up{1} \times (\up{1} \times \ep{0}_{i}) \, , & \text{ for } i = 1,2,3 \, ,
\end{align*}
for some $\up{1}=\up{1}(t)$, $\up{2}=\up{2}(t) \in \R^3$ (see appendix~\ref{app-exp SO3}).

Substituting the expansion\re{smallpsi exp es} in\re{sys-smallpsi}, we find at zeroth order
\begin{equation*}
	\varsigma \, \dot \tau^{[0]} \, \betabody_{0} \times \ep{0}_{i} = 0 \, ,
\end{equation*}
implying that $\taup{0}$ is constant, and that
\begin{equation}
\label{smallpsi-zeroth-sol}
\begin{aligned}
	\ebody_{1} &= \cos \taup{0} \, \betabody_{1} + \varsigma \sin \taup{0} \, \betabody_{2} 
	+ \eps \, (-\sin \taup{1} \, \betabody_{1} + \varsigma \cos \taup{1} \, \betabody_{2}) + \mathcal O (\eps^2) \, , \\
	\ebody_{2} &= -\sin \taup{0} \, \betabody_{1} + \varsigma \cos \taup{0} \, \betabody_{2} 
	- \eps \, (\cos \taup{1} \, \betabody_{1} + \varsigma \sin \taup{1} \, \betabody_{2}) + \mathcal O (\eps^2) \, .
\end{aligned}
\end{equation}

\subsubsection{First order}
Substituting\re{smallpsi-zeroth-sol} in\re{sys-smallpsi}, we obtain at first order
\begin{equation}
	\label{smallpsi-first}
	\varsigma \, \dot \tau^{[1]} \, \betabody_{0} + \dot{\ubody}^{[1]} = - \varsigma \, \Pbody \, (\up{1} \times \betabody_{0}) - \Pbody \, (\cos (a \, t + \taup{0}) \, \betabody_{1} + \varsigma \sin (a \, t + \taup{0}) \, \betabody_{2}) \, .
\end{equation}
We assume that $\up{1} = u_{1} \, \betabody_{1} + u_{2} \, \betabody_{2}$ and solve for $u_{1}$, $u_{2}$. A projection of\re{smallpsi-first} on $\betabody_{1}$ and $\betabody_{2}$ yields
\begin{equation}
	\label{smallpsi-u1}
	\begin{bmatrix}
		u_{1} \\ u_{2}
	\end{bmatrix}
	=
	- a \, (a^2 \, I + A^2)^{-1} \, A
	\begin{bmatrix}
		\varsigma \, \cos (a \, t + \taup{0}) \\
		\sin \, (a \, t + \taup{0})
	\end{bmatrix}
	- (a^2 \, I + A^2)^{-1} \, A^2
	\begin{bmatrix}
		-\varsigma \, \sin(a \, t + \taup{0}) \\ \cos(a \, t + \taup{0})
	\end{bmatrix} \, ,
\end{equation}
where
\begin{equation*}
	A = \begin{bmatrix}
		- \betabody_{1} \cdot \Pbody \betabody_{2} & \betabody_{1} \cdot \Pbody \betabody_{1} \\
		- \betabody_{2} \cdot \Pbody \betabody_{2} & \betabody_{2} \cdot \Pbody \betabody_{1}
	\end{bmatrix} \, .
\end{equation*}
Projecting on $\betabody_{0}$ and substituting\re{smallpsi-u1} therein yields
\begin{equation*}
	\taup{1} = \tilde \tau_{1} \, \cos (a \, t + \taup{0}) + \tilde \tau_{2} \, \sin (a \, t + \taup{0}) \, ,
\end{equation*}
where
\begin{align*}
	\tilde \tau_{1} =& \frac{\varsigma}{a} \,
	\begin{bmatrix}
	\Pbody \, \betabody_{2} \, \cdot \betabody_{0} \\
	-\Pbody \, \betabody_{1} \cdot \betabody_{0}
	\end{bmatrix}
	\cdot \left( - a \, (a^2 \, I + A^2)^{-1} \, A \, \begin{bsmallmatrix} 0 \\ 1 \end{bsmallmatrix} + \varsigma \, (a^2 \, I + A^2)^{-1} \, A^2 \, \begin{bsmallmatrix} 1 \\ 0 \end{bsmallmatrix} \right) \\
	&- \frac{1}{a} \, \Pbody \, \betabody_{2} \cdot \betabody_{0} \\
	\tilde \tau_{2} =& \frac{\varsigma}{a} \,
	\begin{bmatrix}
	\Pbody \, \betabody_{2} \, \cdot \betabody_{0} \\
	-\Pbody \, \betabody_{1} \cdot \betabody_{0}
	\end{bmatrix}
	\cdot \left( \varsigma \, a \, (a^2 \, I + A^2)^{-1} \, A \, \begin{bsmallmatrix} 1 \\ 0 \end{bsmallmatrix} + (a^2 \, I + A^2)^{-1} \, A^2 \, \begin{bsmallmatrix} 0 \\ 1 \end{bsmallmatrix} \right) \\
	&+ \frac{\varsigma}{a} \, \Pbody \, \betabody_{1} \cdot \betabody_{0} \, .
\end{align*}

Note in particular that
\begin{align}
	&\| \up{1} \| \underset{a \to \infty}{\sim} \frac{1}{a} \, , \quad \text{ and} \label{smallpsi-a small}\\
	&\| \up{1} \| \underset{a \to 0}{\sim} 1 \, , \label{smallpsi-a large}
\end{align}
which is consistent with $\ebody_{3} \underset{a \to \infty}{\to} \betabody_{0}$ and $\ebody_{3} \underset{a \to 0}{\to} \begin{bsmallmatrix} \sin \psi \\ 0 \\ \cos \psi \end{bsmallmatrix}$.

\subsubsection{Second order}
At second order, we find that the dynamics is given by
\begin{align*}
	&\varsigma \, \dot \tau^{[2]} \, \betabody_{0} + \dot{\ubody}^{[2]} + \varsigma \, \taup{1} \, \up{1} \times \betabody_{0} + \frac{1}{2} \, \up{1} \times \dot{\ubody}^{[1]} \\
	=& - \varsigma \, \Pbody \, \left( \up{2} \times \betabody_{0} + \frac{1}{2} \, \up{1} \times (\up{1} \times \betabody_{0}) \right) \\
	&- \Pbody \, \left( \cos (a \, t + \taup{0}) \, \up{1} \times \betabody_{1} + \varsigma \sin (a \, t + \taup{0}) \, \up{1} \times \betabody_{2} \right) \, .
\end{align*}
Again we assume that $\up{2} \bot \betabody_{0}$ and we find that $\up{2}$ is an affine combination of $\cos (2 \, a \, t + 2 \, \taup{0})$ and $\sin (2 \, a \, t + 2 \, \taup{0})$ with coefficients depending on $a$.

We use this second order solution only to compare with the numerics.

\subsection{Dynamics of the Magnetic Moment}

The position of the magnetic moment in the magnetic frame is given by
\begin{equation*}
	Q^T \, \mbody = R_{3}^T (a \, t)
	\begin{bmatrix}
		\ebody_{1} \cdot \betabody_{0} \\
		\ebody_{2} \cdot \betabody_{0} \\
		\ebody_{3} \cdot \betabody_{0} \\
	\end{bmatrix} \, .
\end{equation*}
Substituting the $\ebody_{i}$s by their asymptotic expansion for $\sin \psi = \eps \ll 1$, we find that
\begin{equation}
	\label{smallpsi-m mag}
	Q^T \, \mbody
	= R_{3}^T (a \, t + \taup{0})
	\left( 
	\begin{bmatrix} 0 \\ 0 \\ \varsigma \end{bmatrix}
	 + \eps \begin{bmatrix} - u_{2} \\ \varsigma \, u_{1} \\ 0 \end{bmatrix}
	 + \eps^2 \begin{bmatrix} - \up{2} \cdot \betabody_{2} \\ \varsigma \, \up{2} \cdot \betabody_{1} \\ - \varsigma \, \frac{u_{1}^2 + u_{2}^2}{2} \end{bmatrix} \right) + \mathcal O (\eps^3) \, .
\end{equation}

When $1/a \leq \eps$, we use\re{smallpsi-a small} in\re{smallpsi-m mag} to obtain
\begin{equation*}
	Q^T \, \mbody = \begin{bsmallmatrix} 0 \\ 0 \\ \varsigma \end{bsmallmatrix} + \mathcal O (\eps^2)
\end{equation*}
in agreement with our findings for $a \gg 1$ in section~\ref{sec-higha}. When $a \le \eps$, we use\re{smallpsi-a large} to obtain
\begin{equation*}
	Q^T \, \mbody = \begin{bsmallmatrix} 0 \\ 0 \\ \varsigma \end{bsmallmatrix} + \begin{bsmallmatrix} \eps \\ 0 \\ 0 \end{bsmallmatrix} + \mathcal O (\eps^2)
\end{equation*}
in agreement with the prediction for $a \ll 1$ in section~\ref{sec-lowa} since in the magnetic frame, the magnetic field is $[ \eps , 0 , \varsigma \, (1 + \eps^2/2) + \mathcal O (\eps^3) ]^T$.

Moreover, the trajectory of $Q^T \, \mbody$ at first order is a circle with radius $r$ and centre $\mlab_{0}$ satisfying
\begin{align*}
	r &= \frac{\eps \, a}{2 \, \det(a^2 \, I+A^2)} \, \sqrt{c_{0} + c_{1} \, a + c_{2} \, a^2 + c_{3} \, a^3 + c_{4} \, a^4}  \, ,\\
	\mlab_{0} &=
	\begin{bmatrix} 0 \\ 0 \\ \varsigma \end{bmatrix} + \frac{\eps}{2 \, \det(a^2 \, I+A^2)}
	\begin{bmatrix}
		2 \, (\sigma_{1} \, \sigma_{2} \, \cos \iota)^2- a^2 \, (\betabody_{1} \cdot \Pbody \, \Pbody^T \, \betabody_{1} + \betabody_{2} \cdot \Pbody \, \Pbody^T \, \betabody_{2}) \\
		a \, (a^2 + \sigma_{1} \, \sigma_{2} \, \cos \iota) \, (\betabody_{2} \cdot \Pbody \, \betabody_{1} - \betabody_{1} \cdot \Pbody \, \betabody_{2}) \\
		0
	\end{bmatrix} \, ,
\end{align*}
where
\begin{align*}
	c_{0} =& \sigma_{1}^2 \, \sigma_{2}^2 \, \cos^2 \iota\, \left( (\Pbody \betabody_{1} \cdot \betabody_{1} -\Pbody \, \betabody_{2} \cdot \betabody_{2})^2 +(\Pbody \, \betabody_{1} \cdot \betabody_{2} + \Pbody \, \betabody_{2} \cdot \betabody_{1})^2 \right) \, ,\\
	c_{1} =& 2 \, \sigma_{1} \, \sigma_{2} \, \cos \iota \, \big( - \varsigma \, (\Pbody \, \betabody_{1} \cdot \betabody_{1} - \Pbody \, \betabody_{2} \cdot \betabody_{2}) \, (-\Pbody \, \betabody_1 \cdot \betabody_0^2 - \Pbody \, \betabody_2 \cdot \betabody_0^2 + \sigma_{1}^2 + \sigma_{2}^2) \\
	&- 4 \, (\Pbody \, \betabody_1 \cdot \betabody_2 + \Pbody \, \betabody_2 \cdot \betabody_1) \,(\Pbody \, \betabody_1 \cdot \betabody_0) \, (\Pbody \, \betabody_2 \cdot \betabody_0) \big) \, ,\\
c_2 =& -2 \, \sigma_1 \, \sigma_2 \,\cos \iota\, \big((\Pbody\, \betabody_1 \cdot \betabody_1-
\Pbody \, \betabody_2 \cdot \betabody_2)^2 +(\Pbody \, \betabody_1\cdot \betabody_2+\Pbody\,
\betabody_2\cdot \betabody_1)^2\big)\\
	&+ (-\Pbody \, \betabody_1 \cdot \betabody_0^2-\Pbody \, \betabody_2 \cdot \betabody_0^2+\sigma_1^2+\sigma_2^2)^2 + 4 \, (\Pbody \, \betabody_1 \cdot \betabody_0)^2 \, (\Pbody \, \betabody_2 \cdot \betabody_0)^2 \, ,\\
	c_{3} =& 2 \, \varsigma \, (-\Pbody \, \betabody_1 \cdot \betabody_0^2-\Pbody \, \betabody_2 \cdot \betabody_0^2+\sigma_1^2+\sigma_2^2) \, (\Pbody \, \betabody_1 \cdot \betabody_1-\Pbody \, \betabody_2 \cdot \betabody_2) \\
	&+ 4 \, (\Pbody \, \betabody_1 \cdot \betabody_0) \, (\Pbody \, \betabody_2 \cdot \betabody_0) \, (\Pbody \, \betabody_1 \cdot \betabody_2+\Pbody \, \betabody_2 \cdot \betabody_1) \, ,\\
	c_{4} =& (\Pbody \, \betabody_1 \cdot \betabody_1-\Pbody \, \betabody_2 \cdot \betabody_2)^2 + (\Pbody \, \betabody_1 \cdot \betabody_2+\Pbody \, \betabody_2 \cdot \betabody_1)^2 \, ,
\end{align*}
and $\det(a^2 \, I + A^2) = (\sigma_{1} \, \sigma_{2} \, \cos \iota)^2 + a^2 \, \mathrm{tr}(A^2) + a^4$.

This behaviour at small $\sin \psi$ matches the limits found at low and large Mason numbers $a$ in sections~\ref{sec-lowa} and~\ref{sec-higha}. These findings are supported by numerical experiments and can also be observed at larger values of $\sin \psi$ (cf figures~\ref{fig: small psi comp} and~\ref{fig: psi=pi/4}).

\section{Numerical Integration} \label{sec-dni}

We compare analytical predictions of sections~\ref{sec-lowa}-\ref{sec-smallpsi} with solutions obtained by direct numerical integration with integrator ode45 in MATLAB. Using the MATLAB package MatCont~\cite{Dhooge2008}, we also apply numerical continuation starting from known solutions to investigate the existence of periodic solutions to\re{sysQ} across the whole parameter plane.

\subsection{Adapting the system for numerical treatment using quaternions}

We have already cast our system in the three equivalent forms\re{syses},\re{syseB}, and\re{sysQ}. Numerically, we found it easier to work with\re{sysQ} as an ode on SO(3). However, the Lie group SO(3) of rotation matrices has dimension three as opposed to the nine components of a 3 by 3 matrix. Allowing numerical integrators to treat\re{sysQ} efficiently therefore requires using a parametrisation of SO(3). Unit quaternions provide a particularly adapted parametrisation as they both elegantly describe rotations~\cite{Altmann2005} and can be easily processed by numerical integrators as vectors in $\R^4$. In this parametrisation, the ode\re{sysQ} becomes~\cite[]{Dichmann1996}
\begin{equation}
	\label{quaternion form}
	\begin{aligned}
		\dot{q} &= \frac{1}{2} \, F^T(q) \, \ubody (q; a, \psi)
		\, , &
		q \left( 0 \right) &= q_{0} \, ,
	\end{aligned}
\end{equation}
where
\begin{equation*}
	F\left( q \right) =
	\begin{bmatrix}
		\phantom{-}q_{4} & -q_{3} & \phantom{-}q_{2} & -q_{1} \\
		\phantom{-}q_{3} & \phantom{-}q_{4} & -q_{1} & -q_{2} \\
		-q_{2} & \phantom{-}q_{1} & \phantom{-}q_{4} & -q_{3}
	\end{bmatrix}
\end{equation*}
and
\begin{equation*}
	\ubody \left( q; a, \psi \right) = a \, Q(q) \, \begin{bsmallmatrix}0 \\ 0 \\ 1\end{bsmallmatrix} - \Pbody \, Q(q) \, \begin{bsmallmatrix}\sin \psi \\ 0 \\ \cos \psi\end{bsmallmatrix} \, ,
\end{equation*}
with
\begin{equation*}
	Q \left( q \right) =\frac{1}{|q|^2}
	\begin{bmatrix}
		q_{1}^2 -q_{2}^2 -q_{3}^2+q_{4}^2 & 2\left( q_{1} \, q_{2}-q_{3} \, q_{4} \right) & 2\left( q_{1} \, q_{3}+q_{2} \, q_{4} \right) \\
		2\left( q_{1} \, q_{2} + q_{3} \, q_{4} \right) & -q_{1}^2+q_{2}^2-q_{3}^2+q_{4}^2 & 2\left( q_{2} q_{3}-q_{1} \, q_{4} \right) \\
		2\left( q_{1} \, q_{3}-q_{2} \, q_{4} \right) & 2\left( q_{2} \, q_{3}+q_{1} \, q_{4} \right) & -q_{1}^2 - q_{2}^2 + q_{3}^2 + q_{4}^2
	\end{bmatrix} \, .
\end{equation*}
This parametrisation of rotations by quaternion is independent of the norm of the quaternion, and in fact we use only unit quaternions. Solutions of\re{quaternion form} are analytically guaranteed to preserve the norm of the initial condition, which is convenient. However, this independence on the norm also causes the Jacobian of the RHS of\re{quaternion form} to be singular. As this is an issue for numerical continuation, instead of\re{quaternion form} we integrate the modified system
\begin{equation}
	\label{quaternion form corrected}
	\begin{aligned}
		\dot{q} &= \frac{1}{2} \, F^T \left( q \right) \, \ubody \left( q; a, \psi \right) - \frac{1}{2} \, \left( \left| q \right|^2 - 1 \right) \, q
		\, , &
		q\left( 0 \right) &= q_{0} \, .
	\end{aligned}
\end{equation}
Solutions of\re{quaternion form corrected} with unit initial conditions are guaranteed to remain of norm one. They are also solutions of\re{quaternion form}, and the stability of steady states and periodic orbits is kept unchanged.

\subsection{Numerical Solutions}\label{subsec-numsol}

\begin{figure}[h]
	\centering
	\includegraphics[width=275pt]{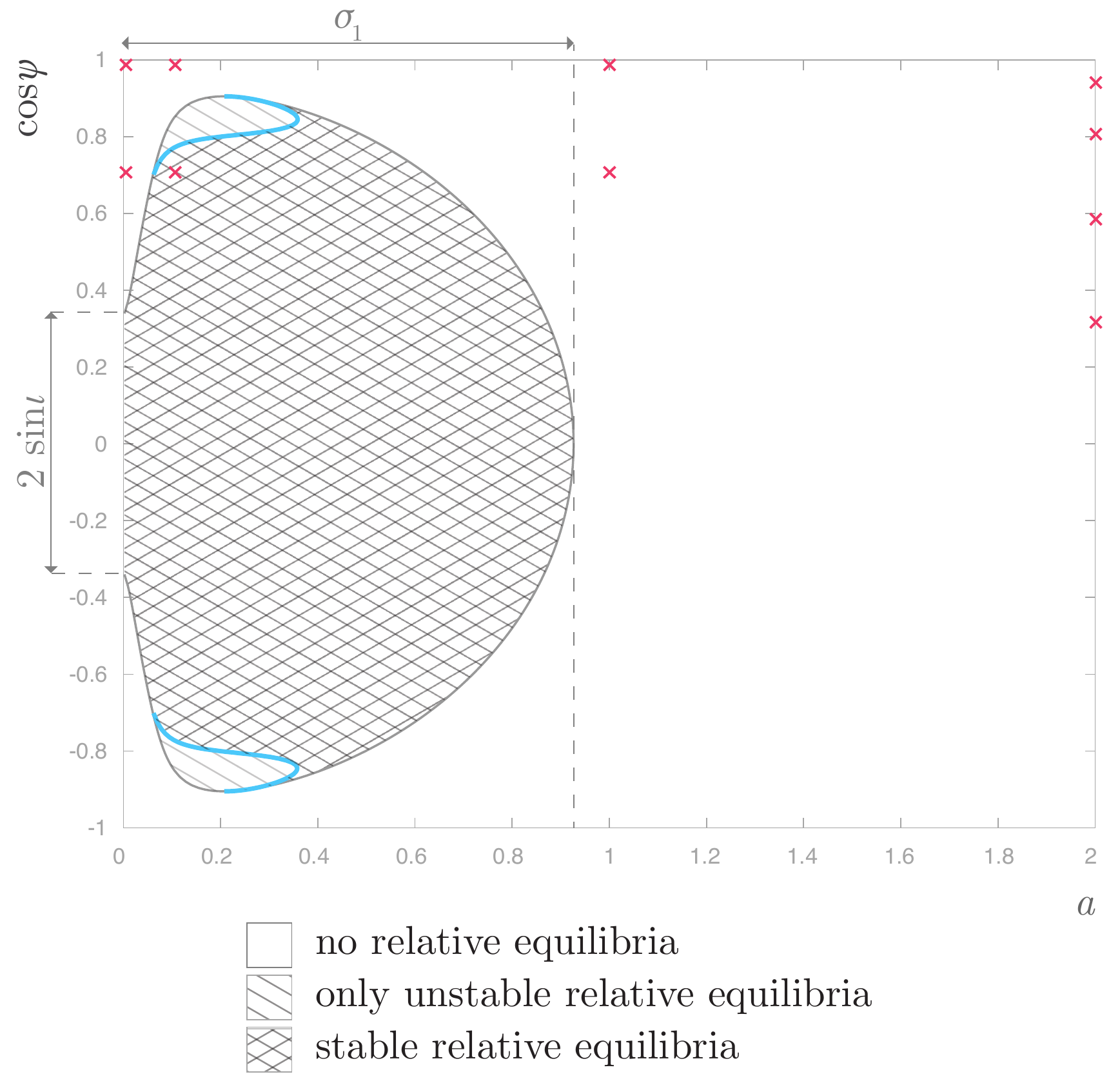}
	\caption{The parameter plane, with regions hatched according to the existence of stable relative equilibria as investigated in~\cite{Ruegg-Reymond2018,Ruegg-Reymond2019} for a specific swimmer (cf appendix~\ref{app-swimmer}). All the parameters for which steady states exist are shown in this figure. We establish the existence of stable periodic solutions to\re{syseB} or equivalently\re{sysQ} for all parameters in the white and single-hatched areas, and for some parameters in the cross-hatched area using numerical methods. The red crosses indicate positions of parameter values used in figures~\ref{fig: a=2}, \ref{fig: small psi comp}, and~\ref{fig: psi=pi/4}.}
	\label{fig: param regimes}
\end{figure}

Direct numerical integration of\re{quaternion form corrected} was performed in MATLAB using the standard integrator ode45~\cite{Shampine1997}. The observed solutions converged either towards stable equilibria or towards stable periodic solutions depending on parameters $a$ and $\psi$, and on initial conditions. Equilibria of\re{quaternion form corrected} are entirely classified and are the object of a separate publication~\cite{Ruegg-Reymond2018}. To explore the existence of periodic solutions for parameters $a$ and $\psi$ for which there are no stable equilibrium, we used numerical continuation (MatCont~\cite{Dhooge2008}) starting from periodic solutions discovered by direct numerical integration. This procedure allowed to find connected sets of periodic orbits seemingly covering the entire parameter plane except for part of the region where stable steady states exist.

We choose here to present the example of a specific swimmer that has the shape of a helical rod. Unless stated otherwise, the features observed for this swimmer are persistent across a range of different helical swimmers, and we conjecture that some of them, in particular the existence of either stable steady states or stable periodic orbits for any pair of parameters $a$ and $\psi$, remain for swimmers of any shape. The study of precisely how solutions depend on the swimmer's shape focusing on helical swimmers will be the subject of a separate publication.

For the example considered (cf appendix~\ref{app-swimmer}), we actually found two distinct sets of periodic orbits, each seemingly covering the entire region of parameter plane where there are no stable equilibria, and part of the region where stable equilibria exist (cf fig.~\ref{fig: param regimes}). They are related to each other through the symmetry of system\re{sysQ} under transformation $Q(t) \mapsto Q(-t) \, R_{2}(\pi)$. These sets intersect along a line seemingly close to $\psi = \pi/2$ corresponding to a stability exchange. One of them contains stable periodic solutions only on one side of the intersection line, and the other one only on the other side. Several other bifurcations corresponding to stability exchanges occur on these sets, notably fold bifurcations of periodic orbits resulting in regions where two distinct stable periodic solutions coexist. Loss of stability also occurs along some branches of periodic solutions as they approach a region with stable steady states, although the type of bifurcation could not be identified with certitude. It is noteworthy that not all branches lose stability as they enter this region.

Other families of periodic orbits that are disconnected from the two sets mentioned above also exist. For instance, there are families of periodic orbits bifurcating from Hopf bifurcations that are part of the set of equilibria. For the swimmer considered here, these families contain stable periodic orbits but this is not the case for other helical swimmers. We cannot rule out the possibility that there exist families of periodic that are neither part of the sets spanning most of the parameter space nor part of the periodic orbits bifurcating from Hopf bifurcations. However we are confident that for the specific swimmer considered here, all the stable periodic solutions were found. Indeed a numerical campaign was set up to investigate the existence of other stable solutions: direct numerical integration was performed for various values of the parameters $a$ and $\psi$ and for 100 randomly chosen initial condition each time. Only solutions corresponding to those already discussed here were found.

\subsection{Comparisons between analytical predictions and numerical solutions}

\begin{figure}[h]
	\centering
	\includegraphics[width=344pt]{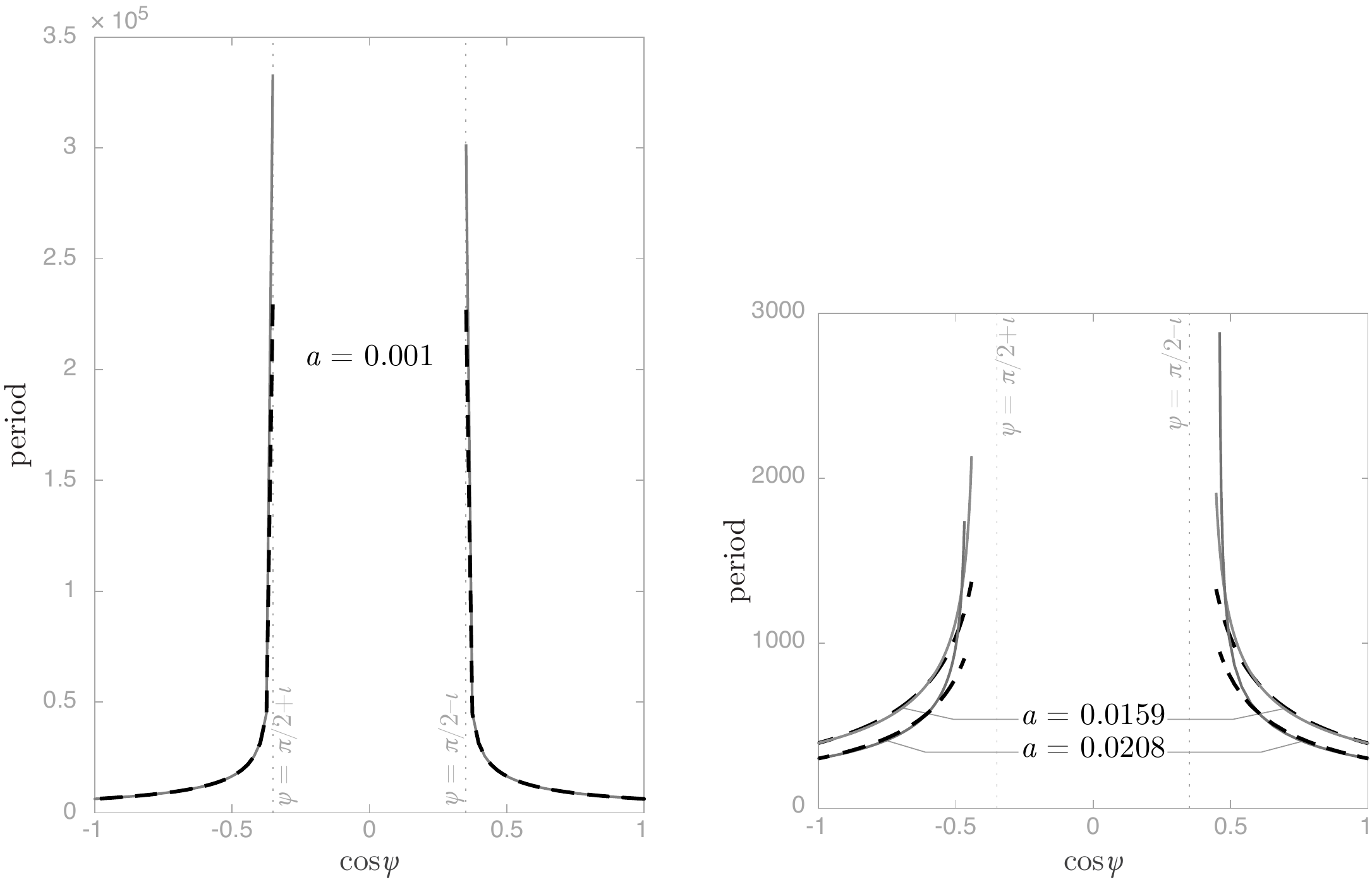}
	\caption{Periods of the periodic solutions vs $\cos \psi$ for several values of $a$\\
	\textcolor{gray1}{\rule{1cm}{1.2pt}} as obtained by numerical computation\\
	\textcolor{gray2}{\hdashrule{1cm}{1.7pt}{2.5mm 2pt}} as predicted analytically in section~\ref{sec-lowa}.\\
	On the left, $a = 10^{-3}$, and the computed period was recovered from solutions obtained by direct numerical integration of\re{quaternion form corrected} for various values of $\psi$. On the right, the computed period was obtained by numerical continuation letting $\psi$ vary starting from a periodic solution found by direct integration.
	}
	\label{fig: small a comp}
\end{figure}

For small Mason number $a$, we obtain limit cycles both numerically and analytically for $\psi \in [0,\pi/2-\iota) \cup (\pi/2+\iota,\pi]$. Figure~\ref{fig: small a comp} displays the periods of analytical and numerical limit cycles for several fixed values of $a$ and varying $\psi\in[0,\pi/2-\iota) \cup (\pi/2+\iota,\pi]$. For $a = 10^{-3}$, the numerical solutions were obtained by direct integration for different values of $\psi$ whereas for $a = 0.0159$ and $a = 0.0208$, the solutions were obtained by numerical continuation letting $\psi$ vary. Note that to numerically characterise the periods, we must take into account that the quaternion parametrisation of SO(3) is a two to one covering: $-q$ and $q$ parametrise the same rotation. Thus there can be symmetric limit cycles~\cite[p. 282]{Kuznetsov2004} in the quaternion coordinates that actually correspond to limit cycles in SO(3) of half the period. All the limit cycles found numerically for small $a$ fall into this category. Therefore we compare the period obtained analytically in\re{periodlowa} with half the period obtained numerically. 

For $a = 10^{-3}$ the relative error between the two is smaller than\footnote{The relative error has finite local maxima, and explodes as it approaches $\pi/2-\iota$ from the left or $\pi/2 + \iota$ from the right. The bounds are chosen taking into account the largest of these local maxima.} $7.6922\cdot10^{-5}$ for $\psi < \pi/2-\iota - 0.1988$ and $\psi > \pi/2 + \iota + 0.2235$. For $\psi$ closer to $\pi/2 \pm \iota$, the relative error  becomes large (the maximal computed value is $0.3688$). This behaviour is expected when trying to approximate a vertical asymptote. 

Note that $a = 0.0159$ and $a = 0.0208$ don't fall into the category of asymptotically small $a$ for this problem. Indeed, the two singular values $\sigma_{1}$ and $\sigma_{2}$ provide characteristic dimensions, and for $a$ to be considered asymptotically small it must verify $a \ll \min \{\sigma_{1}, \sigma_{2} \}$. Here the minimal singular value is $\sigma_{2} = 0.0497$. For $a = 0.0159$ the relative error is smaller than $0.0198$ for $\psi \in [0,\pi/2-\iota-0.2267) \cup (\pi/2+\iota+0.1002,\pi]$. For $a=0.0208$, the relative error is smaller than $0.0291$ for $\psi \in [0,\pi/2-\iota-0.2126) \cup (\pi/2+\iota+0.1297,\pi]$.  That the analysis and numerics fit so well for these values of $a$ illustrates the robustness of the features captured by the asymptotic expansion described in section~\ref{sec-lowa}.

\begin{figure}[h]
	\centering
	\includegraphics[width=275pt]{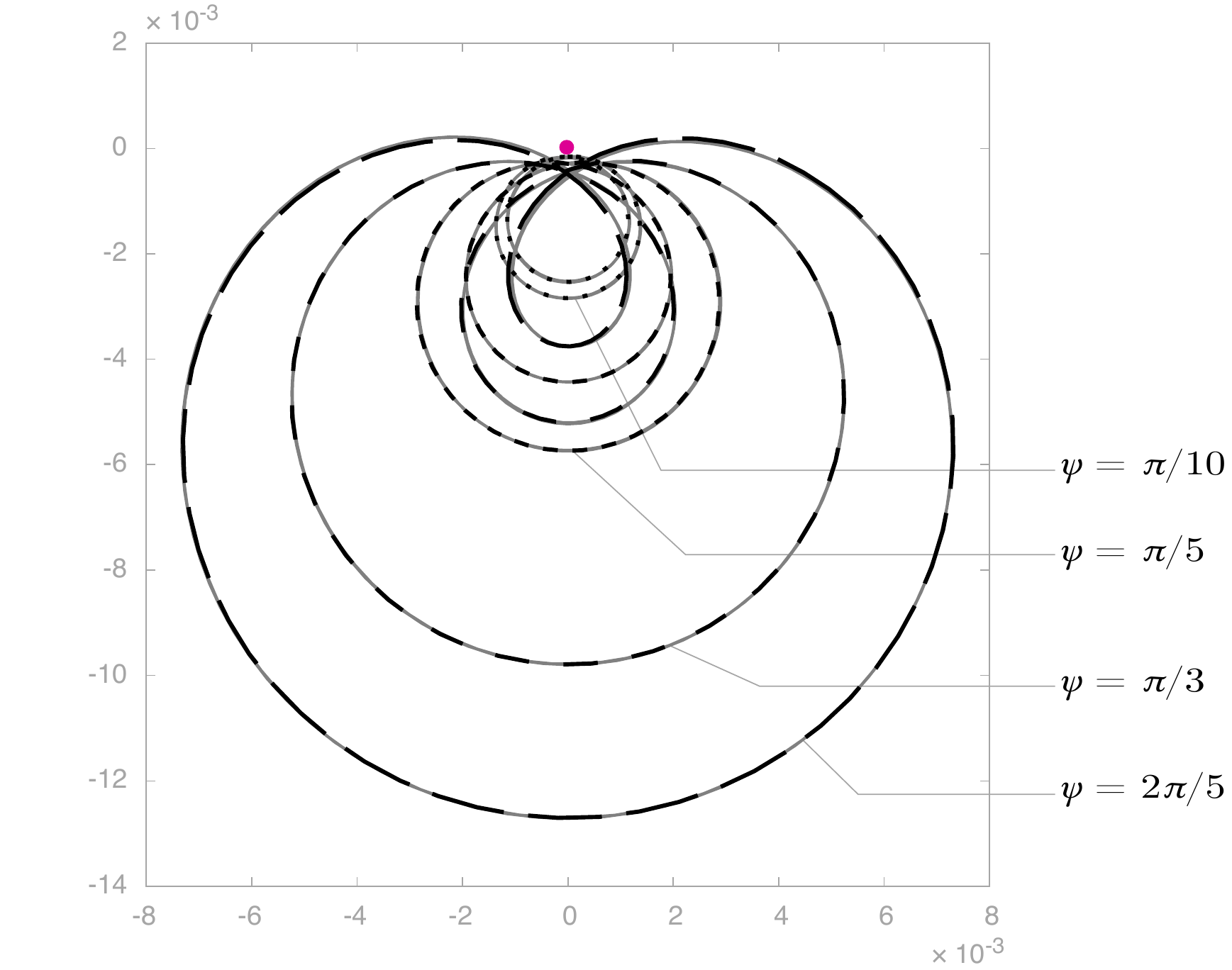}
	\caption{Trajectories of the magnetic moment in the magnetic frame for $a = 100$ and several values of parameter $\psi$\\
	\textcolor{gray1}{\rule{1cm}{1.2pt}} as obtained by numerical integration of\re{quaternion form corrected}\\
	\textcolor{gray2}{\hdashrule{1cm}{1.7pt}{2.5mm 2pt}} as predicted analytically in section~\ref{sec-higha}\\
	\textcolor{magFrame}{$\bullet$} indicates the position of the axis of rotation of the magnetic field $\ebody_{3}$.\\
	This view is a projection perpendicular to $\ebody_{3}$ (notice the scales).
	}
	\label{fig: large a comp}
\end{figure}

\begin{figure}[h]
	\centering
	\includegraphics[width=275pt]{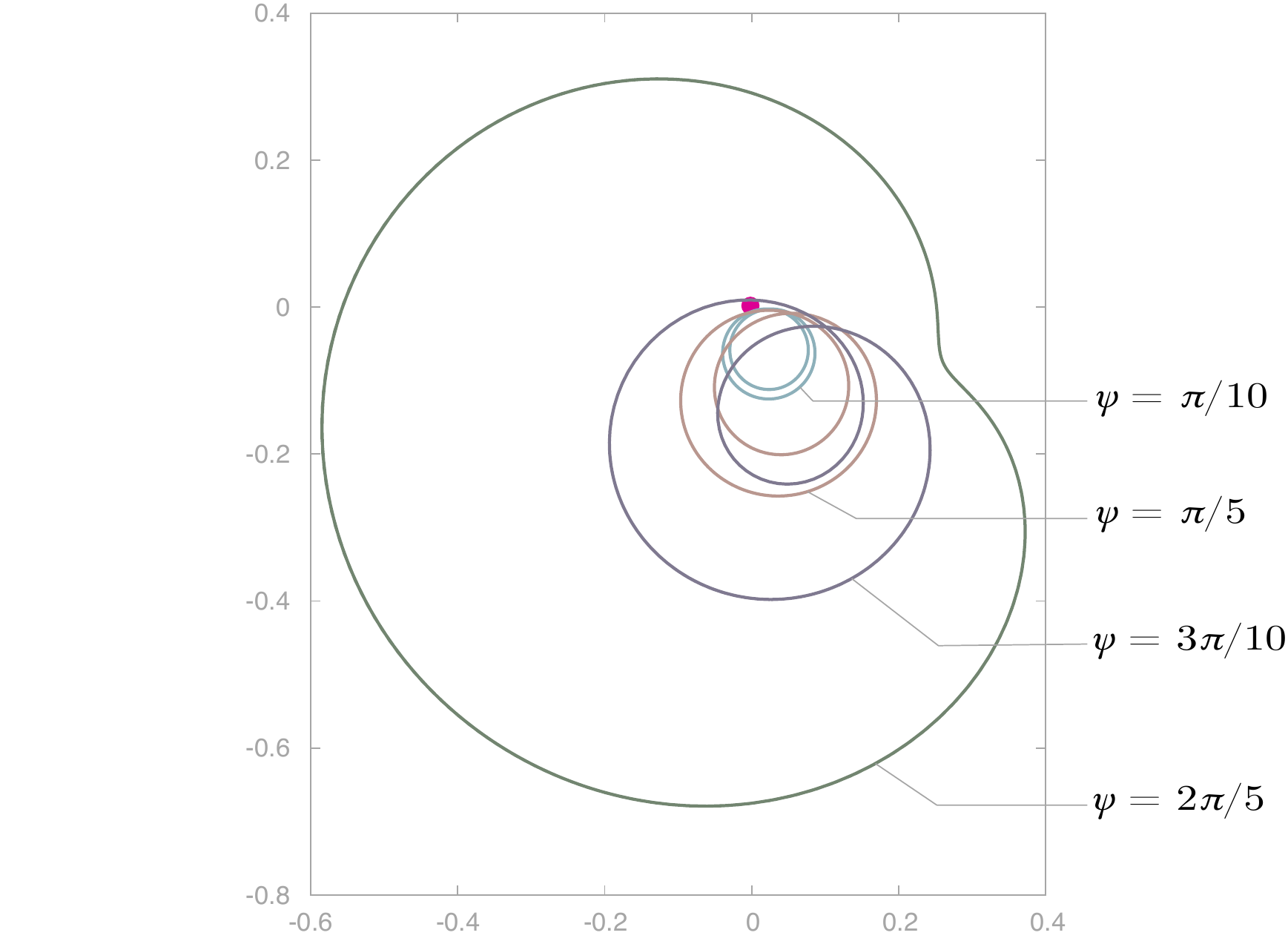}
	\caption{Trajectories of the magnetic moment in the magnetic frame for $a = 2$ and several values of parameter $\psi$ as obtained by numerical integration of\re{quaternion form corrected}\\
	\textcolor{magFrame}{$\bullet$} indicates the position of the axis of rotation of the magnetic field $\ebody_{3}$.\\
	This view is a projection perpendicular to $\ebody_{3}$.\\
	This value of $a$ is out of the scope of the expansion for large $a$ in section~\ref{sec-higha}, but the behaviour we observe is qualitatively similar as the magnetic moment rotates around the average magnetic field. Note that the approximation is better for lower values of $\psi$.
	}
	\label{fig: a=2}
\end{figure}

In order to visualise the agreement between solutions obtained numerically and analytically in both the large $a$ and small $\sin\psi$ regimes, we find that it is useful to view the curve described by the magnetic moment $\mbody$ in the magnetic frame, that is in the frame locked to the rotating magnetic field $\Bbody$. In the large $a$ regime, figure~\ref{fig: large a comp} exhibits a remarkable agreement between the curves obtained as a first order expansion as in section~\ref{sec-higha} and by direct numerical integration for $a = 100$ and different values of $\psi$. In figure~\ref{fig: a=2} we show the results of direct numerical integrations for $a = 2$ and different values of $\psi$. Although $a = 2$ is not large enough to be in the large $a$ regime, we can observe that the magnetic moment stays close to the average magnetic field, especially for small values of $\psi$.

\begin{figure}[h]
	\centering
	\includegraphics[width=275pt]{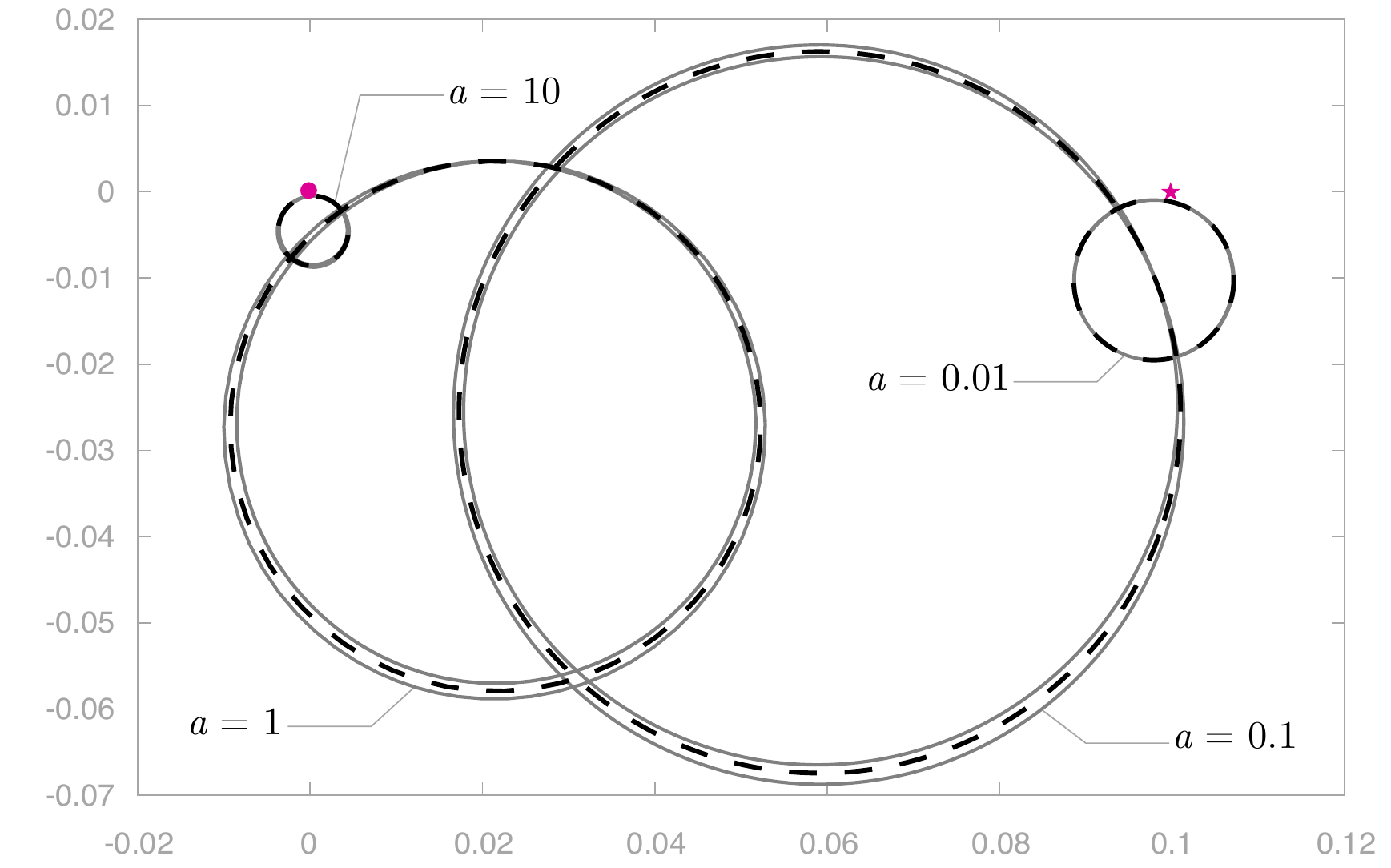}\\
	\vspace{2\baselineskip}
	\includegraphics[width=275pt]{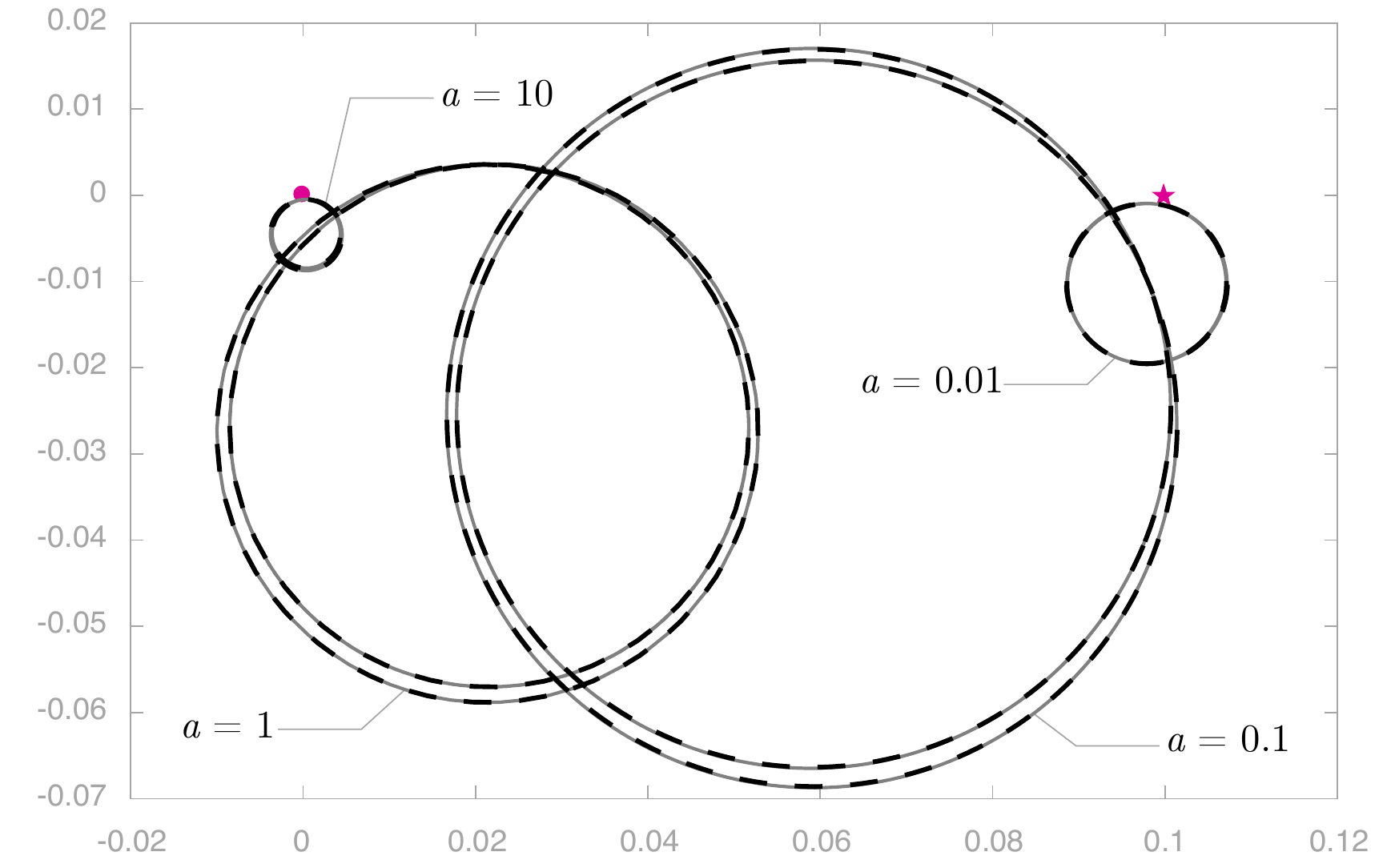}
	\caption{Trajectories of the magnetic moment in the magnetic frame for $\psi = 0.1$ and various Mason numbers $a$\\
	\textcolor{gray1}{\rule{1cm}{1.2pt}} as obtained by numerical integration of\re{quaternion form corrected}\\
	\textcolor{gray2}{\hdashrule{1cm}{1.7pt}{2.5mm 2pt}} as predicted analytically in section~\ref{sec-smallpsi} (top: first order, bottom: second order).\\
	\textcolor{magFrame}{$\bullet$} indicates the position of the axis of rotation of the magnetic field $\ebody_{3}$.\\
	\textcolor{magFrame}{\scriptsize $\bigstar$} indicates the position of the magnetic field $\Bbody$.\\
	This view is a projection perpendicular to $\ebody_{3}$.
	}
	\label{fig: small psi comp}
\end{figure}

\begin{figure}[h]
	\centering
	\includegraphics[width=275pt]{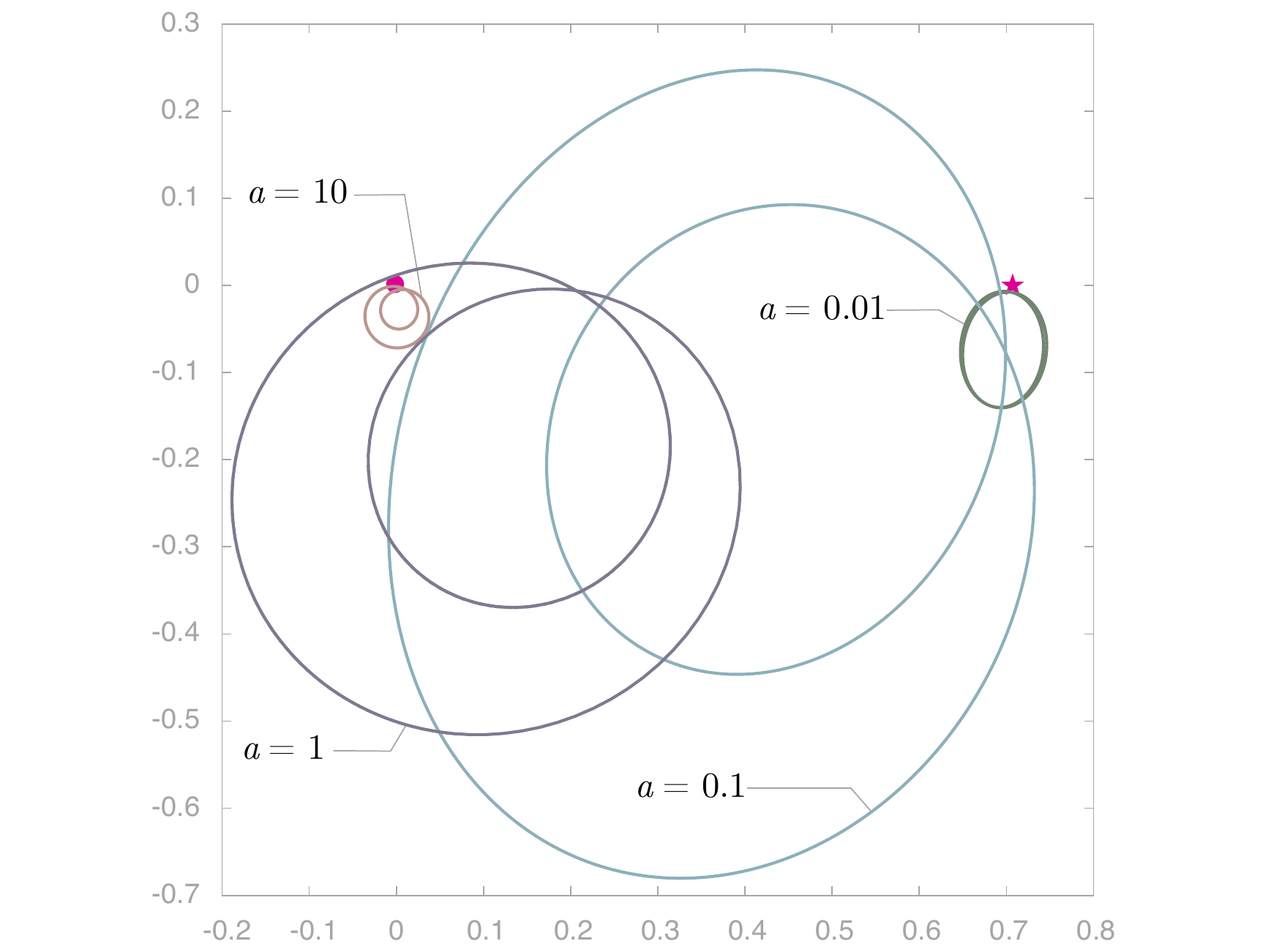}
	\caption{Trajectories described by the magnetic moment in the magnetic frame for $\psi = \pi/4$ and several values of parameter $a$ as obtained by numerical integration of\re{quaternion form corrected}\\
	\textcolor{magFrame}{$\bullet$} indicates the position of the axis of rotation of the magnetic field $\ebody_{3}$.\\
	\textcolor{magFrame}{\scriptsize $\bigstar$} indicates the position of the magnetic field $\Bbody$.\\
	This view is a projection perpendicular to $\ebody_{3}$.\\
	This value of $\psi$ is out of the scope of the expansion for small $\sin \psi$ in section~\ref{sec-smallpsi}, but the behaviour we observe is qualitatively similar.
	}
	\label{fig: psi=pi/4}
\end{figure}

The small $\sin \psi$ regime is displayed in figure~\ref{fig: small psi comp} for $\psi = 0.1$. Note that the circles predicted as a first order approximation give a good estimate of the trajectories, and that the second order prediction virtually overlaps the numerical solutions. For $\psi = \pi/4$, that is outside the scope of this asymptotic expansion, the curves described by $\mbody$ in the magnetic frame have a qualitatively similar behaviour, with a mean position moving from $\Bbody$ to $\ebody_{3}$ as $a$ increases, and a decreasing amplitude in both the low and large $a$ limits (cf fig.~\ref{fig: psi=pi/4}).

Numerical solutions to\re{quaternion form corrected} at large $a$ and at small $\sin\psi$ were found to be periodic. The corresponding analyses in sections~\ref{sec-higha} and~\ref{sec-smallpsi} did not predict periodicity but are consistent with it. Numerical periodic solutions corresponding to the small $a$ and large $a$ regimes of section~\ref{sec-lowa} and~\ref{sec-higha} are connected numerically by branches of periodic orbits that for small $\sin\psi$ correspond to the analytical regime of section~\ref{sec-smallpsi}.

\section{Conclusion}

Using analytical and numerical methods, we studied out-of-equilibrium solutions for the motion of a rigid body in Stokes flow submitted to a steadily rotating external magnetic field uniform in space. The equations governing this dynamics is given in three equivalent forms~(\ref{syses}-\ref{syseB}) depending on two parameters: the Mason number $a$ and conical angle $\psi$. We study analytically their solutions using asymptotic expansions in three different regimes: in the limit $a\to 0$ (section~\ref{sec-lowa}), in the limit $a \to \infty$ (section~\ref{sec-higha}), and in the limit $\sin \psi \to 0$ (section~\ref{sec-smallpsi}).

In the limit of small $a$, which corresponds to either slowly rotating or strong magnetic field, the body aligns its magnetic moment $\mbody$ giving the dipole axis with the magnetic field. We find a limiting angle $\iota$ such that the motion in this regime is in relative equilibrium for $\psi \in (\pi/2-\iota,\pi/2+\iota)$ and in relative periodic motion elsewhere. This periodic motion arises as the body rotates about its magnetic moment. We give the period of this solution explicitly. The change of regime as $\psi$ approaches the limit $\pi/2 \pm \iota$ occurs through a periodic solution of infinite period.

In the limit of large $a$, which corresponds to either rapidly rotating or weak magnetic field, the body aligns its magnetic moment $\mbody$ with the average magnetic field. It exhibits a residual rotation about its magnetic moment. The limit $\sin \psi \to 0$ gives a continuous change between the small and large $a$ regimes, with magnetic moment $\mbody$ shifting from alignment with the magnetic field $\Bbody$ to alignment between the averaged magnetic field as $a$ goes from $0$ to $\infty$.

We used numerical integration and continuation to assess the validity of our analytical results on a particular example. The predicted features are in good qualitative agreement even for parameter values outside the scope of our three asymptotic expansions. Numerical continuation allowed to show that the out-of-equilibrium solutions at large $a$ and small $\sin \psi$ are limit cycles of the dynamics\re{sysQ}. Therefore at all parameter values for the swimmer tested, we find either stable periodic orbits or stable equilibria.

Comparison of our results with experiments is limited by two main factors: we assume here that the swimmers are permanently magnetic, and that they are neutrally buoyant. To our knowledge, these two factors have not arised simultaneously in experiments with magnetic micro-swimmers.

\subsection*{Acknowledgements}

We thank Oscar Gonzalez for providing the code used to compute the drag matrix used in numerical simulations. We also thank Andrew Petruska for the inspiring discussions.

\appendix
\section{Asymptotic expansion on SO(3)}\label{app-exp SO3}

In sections~\ref{sec-higha} and~\ref{sec-smallpsi} we used an asymptotic expansion on SO(3) to obtain\re{nearidentity} and\re{expandc}, and\re{smallpsi exp es} respectively. We show here that an asymptotic expansion on a right-handed frame $\{\ebody_{1}, \ebody_{2}, \ebody_{3} \}$ takes the form
\begin{equation*}
	\ebody_{i} = \ep{0}_{i} + \eps \, \up{1} \times \ep{0}_{i} + \eps^2 \, \left( \up{2} \times \ep{0} + \frac{1}{2} \up{1} \times (\up{1} \times \ep{0}) \right) + \mathcal O (\eps^3) \, .
\end{equation*}
Everything used in this appendix belongs to classical mathematical knowledge, but to the best of our knowledge, the derivation of generic asymptotic expansions in SO(3) can't be found in the literature. We provide it here to ease the progression of the interested reader through sections~\ref{sec-higha} and~\ref{sec-smallpsi}.

We use the matrix form $R = \bma \ebody_{1} & \ebody_{2} & \ebody_{3} \ema \in$~SO(3) of a right-handed frame given by vectors $\ebody_{i}$ for $i = 1,2,3$. Suppose we have a curve on SO(3) given by $\eps \mapsto R(\eps)$. Around $\eps = 0$, this curve is approximated by its Taylor expansion
\begin{equation}
	\label{R approx}
	R(\eps) = \Rex{0} + \eps \, \Rex{1} + \eps^2 \, \Rex{2} + \mathcal O (\eps^3) \, ,
\end{equation}
with
\begin{align}
	\label{def R012}
	\Rex{0} &= R(0) \, , & \Rex{1} &= \left. R' \right|_{\eps = 0} \, , && \text{ and } & \Rex{2} &= \frac{1}{2} \left. R'' \right|_{\eps = 0} \, ,
\end{align}
where $'$ denotes the derivative by $\eps$. Since $R(\eps) \in$~SO(3) for all $\eps$, we have the identity
\begin{equation*}
	R(\eps) \, R^T(\eps) = I \, ,
\end{equation*}
and differentiating it with respect to $\eps$ yields
\begin{equation*}
	R'(\eps) \, R^T(\eps) + R(\eps) \, R^{\prime \, T}(\eps) = 0 \, .
\end{equation*}
This implies that $R' \, R^T$ is a skew-symmetric matrix, so there exists $\ubody = \ubody(\eps) \in \R^3$ such that
\begin{align}
	\label{R ex 1st order}
	\cross[\ubody] &= R' \, R^T && \Leftrightarrow & R' &= \cross[\ubody] R \, .
\end{align}
Equivalently, the three columns $\ebody_{i}$ ($i = 1,2,3$) of $R$ satisfy
\begin{equation*}
	\ebody_{i}' = \ubody \times \ebody_{i} \, .
\end{equation*}
Differentiating again we obtain
\begin{align}
	\label{ex 2nd order}
	R'' &= \cross[\ubody'] \, R + \cross[\ubody]^2 \, R \, && \text{ and } & \ebody_{i}'' &= \ubody' \times \ebody_{i} + \ubody \times (\ubody \times \ebody_{i}) \, .
\end{align}

Expanding $\ubody$ in $\eps$ as $\ubody = \ue{1} + 2 \, \eps \, \ue{2} + \mathcal O (\eps^2)$, and substituting\re{R ex 1st order} and\re{ex 2nd order} in\re{def R012} we find that
\begin{align*}
	\Rex{1} &= \cross[\ue{1}] \, \Rex{0} \, , && \text{ and } & \Rex{2} &= \cross[\ue{2}] \, \Rex{0} + \frac{1}{2} \, \cross[\ue{1}]^2 \, \Rex{0} \, .
\end{align*}
Thus obtain the expansion of $R$ in $\eps$\re{R approx} takes the form
\begin{equation*}
	R = \Rex{0} + \eps \, \cross[\ue{1}] \, \Rex{0} + \eps^2 \, \cross[\ue{2}] \, \Rex{0} + \mathcal O (\eps^3) \, .
\end{equation*}
Equivalently the frame vectors $\ebody_{i}$ ($i = 1,2,3$) satisfy 
\begin{equation*}
	\ebody_{i} = \ep{0}_{i} + \eps \, \up{1} \times \ep{0}_{i} + \eps^2 \, \left( \up{2} \times \ep{0} + \frac{1}{2} \up{1} \times (\up{1} \times \ep{0}) \right) + \mathcal O (\eps^3) \, .
\end{equation*}

\section{Numerical data for example swimmer}\label{app-swimmer}

The swimmer used as an example for this paper is a body with the shape of a helical rod with radius $r = 0.1330$, pitch $p = 1.1076$, total arc-length $L = 4.1628$, and rod radius~0.0936 in non-dimensional units chosen so that its radius of gyration is~1. The magnetic moment is $\mbody =  \left( 0, 0.1736, 0.9848 \right)^T$ in a body frame that is chosen so that the rod's centreline is given by
\begin{equation*}
	s \mapsto
	\begin{bmatrix}
		r \, \cos s \\ r \, \sin s \\ \frac{p}{2 \pi} \, s
	\end{bmatrix}
	\quad \text{for } s \in \left[ 0, \frac{L}{\sqrt{r^2 + (p/2\pi)^2}} \right] \, .
\end{equation*}
The rod is capped at both ends by half spheres. Assuming uniform density, we compute the drag matrix $\dragbody$ by computing the resultant loads due to Stokes flow at the center of mass. To do so, we use a code provided to us by Oscar Gonzalez based on a Nyström method for exterior Stokes flow~\cite{Gonzalez2009,Li2013,Gonzalez2015}. The drag matrix obtained for the example swimmer in the body frame described above is
\begin{equation*}
	\dragbody =
	\begin{bsmallmatrix*}[r]
		12.4654  &  0.0000  &  0.0000  &  0.1433  &  0.0000 &  -0.0000 \\
		-0.0000  & 12.4815  & 0.0582  &  0.0000  &  0.0122  &  0.1178 \\
		-0.0000  &  0.0577   & 9.2808  &  0.0000 &  -0.5607  & -0.2158 \\
		0.1427 &  -0.0000 &  -0.0000  & 20.1070  & -0.0000  &  0.0000 \\
		-0.0000  &  0.0116 &  -0.5610 &  -0.0000  & 20.1725  &  0.4032 \\
		0.0000  &  0.1179  & -0.2158   & 0.0000  &  0.4031  &  1.0196
	\end{bsmallmatrix*}.
\end{equation*}
Since the result is not exactly symmetric, we correct this by using $\Mob = \frac{1}{2} \Drag^{-1} + \frac{1}{2} \Drag^{-T}$ in our computations.

\bibliographystyle{unsrt}
\bibliography{bibliography}

\end{document}